\documentclass[12pt, reqno]{amsart}
\usepackage{amsthm}
\usepackage{amssymb}
\usepackage{amsmath}

\def\1{\hbox{1\kern-.35em\hbox{1}}}

\newtheorem{theorem}{Theorem}[section]
\newtheorem*{theorem*}{Theorem}
\newtheorem{lem}[theorem]{Lemma}

\newtheorem*{proposition*}{Proposition}

\newtheorem{definition}[theorem]{Definition}

\numberwithin{equation}{section}

\newcommand{\bea}{\begin{eqnarray}}
\newcommand{\eea}{\end{eqnarray}}
\newcommand{\be}{\begin{eqnarray*}}
\newcommand{\ee}{\end{eqnarray*}}

\newcommand{\Z}{{\mathbb Z}}

\newcommand{\C}{{\mathbb C}}

\def \pf{\noindent {\bf Proof: \,}}

\def\qed{\hfill\mbox{$\Box$}}

\def\al{\alpha}

\def\dis{\displaystyle}

\def\Z{\mathbb{Z}}

\def\C{\mathbb{C}}

\def\tau{\pi}

\numberwithin{equation}{section}

\begin{document}

\title[Representations for the non-graded Virasoro-like algebra]
{Representations for the non-graded Virasoro-like algebra$^{*}$}
\author[Gao]{Shoulan Gao}
\address{Department of Mathematics, Shanghai Jiaotong University, Shanghai
200240, P.R.China} \email{gaoshoulan@sjtu.edu.cn}

\author[Jiang]{Cuipo Jiang$^{\dag}$}
\address{Department of Mathematics, Shanghai Jiaotong University, Shanghai
200240, P.R.China} \email{cpjiang@sjtu.edu.cn}

\begin{abstract}
It is proved in this paper that  an irreducible module over the
non-graded Virasoro-like algebra $\mathcal{L}$, which satisfies a
natural condition, is a $GHW$ module or uniformly bounded.
Furthermore, we give the complete classification of indecomposable
$\mathcal{L}$-modules $V=\bigoplus_{m,n\in\Z}\mathbb{C}v_{m,n}$
which satisfy $L_{r,s}v_{m,n}\subseteq
\mathbb{C}v_{r+m,s+n+1}+\mathbb{C}v_{r+m,s+n}$.
\end{abstract}

\thanks{$^{*}$ Supported in part by China NSF grant
10571119.  }

\thanks{$^{\dag}$ Corresponding author: cpjiang@sjtu.edu.cn }
\maketitle

\maketitle

%
\section{\bf Introduction}
%

The Virasoro algebra  $Vir$ is the universal central extension of
the Witt algebra with a basis $\{L_{i}, c |\; i\in \mathbb{Z}\}$
such that for all $i, j\in \mathbb{Z}$
$$[L_{i},L_{j}]=(j-i)L_{i+j}+\displaystyle\frac{i^{3}-i}{12}\delta_{i+j,0}c,
\;\;[L_{i},c]=0 .$$ It plays an important role in many areas of
mathematics and  physics. Over the past  decades many authors have
studied the representation of $Vir$ (see for example \cite{CP,K, KS,
MP}). One of the most important results proved by O. Mathieu is that
an irreducible Harish-Chandra module over $Vir$ is a highest weight
module, a lowest weight module or a module of intermediate series
(\cite{M}). And it is known that a module of intermediate series
over $Vir$ is one of
 $A_{a, b},A(a'), B(a')$ or one of their quotients or submodules for
suitable $ a, b\in\mathbb{C},a'\in\mathbb{C}\bigcup\{\infty\}$,
where $A_{a, b},A(a'), B(a')$ all have a basis
$\{v_{j}|j\in\mathbb{Z}\}$, such that $c$ acts trivially and for all
$i,j\in\mathbb{Z}$
\begin{eqnarray*}
A_{a, b}:& &L_{i}v_{j}=(a+bi+j)v_{i+j};
\\
A(a'):& &
L_{i}v_{j}=(i+j)v_{i+j},~~j\neq0,\;\;\;L_{i}v_{0}=i[1+(i+1)a']v_{i},\;
\; a'\in\C;
\\
B(a'):& & L_{i}v_{j}=jv_{i+j},~~j\neq
-i,\;\;\;L_{i}v_{-i}=-i[1+(i+1)a']v_{0},\; \; a'\in\C;
\\
A(\infty):& & L_{i}v_{j}=(i+j)v_{i+j},~~j\neq0,\;\;\;
L_{i}v_{0}=i(i+1)v_{i};
\\
B(\infty):& & L_{i}v_{j}=jv_{i+j},~~j\neq -i,\;\;\;
L_{i}v_{-i}=-i(i+1)v_{0}.
\end{eqnarray*}

Up to now many authors have studied various types of generalizations
of  the Virasoro algebra, such as the Virasoro-like algebra and its
q-analog, generalized Witt algebras, the higher rank Virasoro
algebra and the super-Virasoro algebra and so on, see for example
(\cite{DZh, KPS, JM, LT, PZ, Su1, Su2, Zh}). However, all of these
Lie algebras are graded. Due to the important applications in the
theory of Hamiltonian operator and vertex operator algebras,
infinite dimensional non-graded Lie algebras have been studied
(\cite{ SXZ, SZ2, SZ, X}). In \cite{SZ2}, Su and Zhao calculated the
second cohomology groups of Lie algebras of generalized Witt type
and introduced the non-graded Virasoro-like Lie algebra
$W\widetilde{(\Gamma)}$ ( $\Gamma$ is an additional subgroup of a
field of characteristic zero ) with a basis $\{L'_{\alpha,i},
c\;|\;\alpha\in\Gamma,\; i\in{\mathbb{Z}}\}$ and the product
\begin{eqnarray*}
& [L'_{\alpha,i},L'_{\beta,j}]& =(\beta-\alpha)L'_{\alpha+\beta,i+j}
+(j-i)L'_{\alpha+\beta,i+j-1}
\\
& &+\frac{1}{12}\delta_{\alpha+\beta, 0}[\delta_{i+j,
-1}\alpha^{3}+3i\delta_{i+j, 0}\alpha^{2}+3i(i-1)\delta_{i+j,
1}\alpha+i(i-1)(i-2)\delta_{i+j, 2}]c,
\\
&[L'_{\alpha,\;i},\quad c] &=0.
\end{eqnarray*}

Since these Lie algebras are non-graded and  have no Cartan
subalgebras, we cannot define their weight modules as for  graded
Lie algebras. It makes the representation theory more difficult to
study. So far not much has been achieved related to these non-graded
Lie algebras. In this paper, we study representations of ${\mathcal
L}=W\widetilde{(\Gamma)}$  over $\mathbb{C}$ for $\Gamma=\Z$.

 Let $L_{\alpha,\;\beta}=L'_{\alpha,\;\beta+1}$. Then
 the non-graded Virasoro-like Lie algebra  $\mathcal{L}$ is linearly spanned by $\{L_{\alpha,\beta},
c\;|\;\alpha,\beta\in{\mathbb{Z}}\}$ with the following product
\begin{eqnarray*}
& [L_{\alpha_{1},\beta_{1}},L_{\alpha_{2},\beta_{2}}]&
=(\alpha_{2}-\alpha_{1})L_{\alpha_{1}+\alpha_{2},\beta_{1}+\beta_{2}+1}
+(\beta_{2}-\beta_{1})L_{\alpha_{1}+\alpha_{2},\beta_{1}+\beta_{2}}
\\
& &+\frac{1}{12}\delta_{\alpha_{1}+\alpha_{2},
0}[\delta_{\beta_{1}+\beta_{2}, -3}\alpha_{1}^{3}
+3\delta_{\beta_{1}+\beta_{2}, -2}(\beta_{1}+1)\alpha_{1}^{2}
\\
& &+3\delta_{\beta_{1}+\beta_{2},
-1}\beta_{1}(\beta_{1}+1)\alpha_{1}+\delta_{\beta_{1}+\beta_{2},
0}\beta_{1}(\beta_{1}^{2}-1)]c,
\\
&[L_{\alpha,\;\beta},\quad\quad  c] &=0.
\end{eqnarray*}

Through out the  paper, we study  ${\mathcal L}$-modules
$V=\bigoplus_{r,s\in\Z}V_{r,s}$   such that $\dim V_{r,s}<\infty$
and $L_{a,b}V_{r,s}\subseteq
V_{r+a,s+b+2}+V_{r+a,s+b+1}+V_{r+a,s+b}.$ We prove  in Section 2
that if $V$ is irreducible, then $V$ is either a $GHW$ module or
uniformly bounded. In sections 3-5, we study uniformly bounded
modules of $\mathcal{L}$  and  give the complete classification of
indecomposable ${\mathcal L}$-modules $V=\bigoplus_{m,n\in\Z}V_{m,
n}$ which satisfy the following two conditions: (1) $\dim V_{m,
n}\leq 1$. (2) $L_{\alpha, \beta}V_{m, n}\subseteq V_{\alpha+m,
\beta+n+1}+V_{\alpha+m, \beta+n}$.


\section{\bf $GHW$ Modules}
\label{GHW Modules}
\label{sub2-c-related}

\begin{definition}\label{2.1} Let $\al_{1},\al_{2}\in\Z\times \Z$.
$\{\al_{1},\al_{2}\}$ is called  a $\Z$-basis of $\Z\times\Z$, if
for each $\al\in\Z\times\Z$, $\al=k_{1}\al_{1}+k_{2}\al_{2}$ for
some $k_{1},k_{2}\in\Z$.
\end{definition}
Set
$$e_{1}=(1,0), \ e_{2}=(0,1).$$

In this paper, we always consider the ${\mathcal L}$-modules
$V=\bigoplus_{r,s\in\Z}V_{r,s}$ such that $\dim V_{r,s}<\infty$
and
\begin{equation}\label{ea22.1}L_{a,b}V_{r,s}\subseteq
V_{r+a,s+b+2}+V_{r+a,s+b+1}+V_{r+a,s+b}.
\end{equation}

\begin{definition}\label{2.2} Let $V=\bigoplus_{r,s\in\Z}V_{r,s}$ be
a ${\mathcal L}$-module.  $V$ is called a $GHW$ module if $V$ is
generated by a vector $v$ and there is a $\Z$-basis
$\{\al_{1},\al_{2}\}$ of $\Z\times\Z$ such that
$$L_{\al}v=0,$$
for all
$\al=k_{1}\al_{1}+k_{2}\al_{2}+k(1-\delta_{k_{1}k_{2},0})e_{2}$
with $k_{1},k_{2}\in\Z_{+}$ and $0\leq k\leq k_{1}+k_{2}$.
\end{definition}

\begin{lem}\label{2.3} Let $V=\bigoplus_{r,s\in\Z}V_{r,s}$ be an irreducible
${\mathcal L}$-module and $m,n\in\Z$ be such that $|m|\geq 2$ and
$n\neq 0$. Denote by ${\mathcal S}_{m,n}$ the Lie subalgebra of
${\mathcal L}$ generated by $\{L_{m+i,n+j}|i=1,0,-1, j=3,0,-3\}$.
If there exists $0\neq v\in V$ such that ${\mathcal S}_{m,n}\cdot
v=0$, then $V$ is a $GHW$ ${\mathcal L}$-module.
\end{lem}

\pf {\bf Case 1.} $m,n>0$.  Set
$$\al_{1}=(nm-1)e_{1}+n^{2}e_{2},\ \
\al_{2}=m^{2}e_{1}+(mn+1)e_{2}.$$ Then $\al_{1},\al_{2}$ is a
$\Z$-basis of $\Z\times\Z$. Since
$$[L_{m,n}, L_{m-1,n}]=-L_{2m-1,2n+1},$$
$$
[L_{m,n+3}, L_{m-1,n-3}]=-L_{2m-1,2n+1}-6L_{2m-1,2n},$$ we have
$L_{2m-1,2n}, L_{2m-1,2n+1}\in {\mathcal S}_{m,n}$. By the fact
that
$$
[L_{m,n}, L_{m-1,n-3}]=-L_{2m-1,2n-2}-3L_{2m-1,2n-3},$$
$$
[L_{m-1,n}, L_{m,n-3}]=L_{2m-1,2n-2}-3L_{2m-1,2n-3},$$ we have
$L_{2m-1,2n-2}, L_{2m-1,2n-3}\in {\mathcal S}_{m,n}$. Using the
following Lie bracket relations
$$
[L_{m,n}, L_{m-1,n+3}]=-L_{2m-1,2n+4}-3L_{2m-1,2n+3},$$
$$
[L_{m-1,n}, L_{m,n+3}]=L_{2m-1,2n+4}-3L_{2m-1,2n+3},$$ we deduce
that $L_{2m-1,2n+4}, L_{2m-1,2n+3}\in S_{m,n}$. Therefore,
\begin{equation}\label{ea2.1}\{L_{2m-1,2n+k}|k=-3,-2,0,1,3,4\}\subseteq {\mathcal S}_{m,n}.
\end{equation}
Similarly, we have \begin{equation}\label{ea2.2}\{L_{2m,2n+k},
L_{2m+1,2n+k}|k=-3,-2,0,1,3,4\}\subseteq {\mathcal S}_{m,n}.
\end{equation}
 In
general, for $k\geq 3$,  we have
\begin{equation}\label{ea2.3}\{L_{km+i,kn+j} |i=-1,0,1;-3(k-1)\leq j\leq
4(k-1)\}\subseteq {\mathcal S}_{m,n}. \end{equation} Since $m\geq
2$, by (\ref{ea2.1})-(\ref{ea2.2}) and (\ref{ea2.3}), considering
$k=n$ and $k=m$ respectively, we have
$$
\{L_{mn-1,n^{2}}, L_{mn,n^{2}}, L_{m^{2}+1,mn+1},
L_{m^{2},mn+1}\}\subseteq {\mathcal S}_{m,n}.$$ Then it is easy to
deduce that $L_{\al}\subseteq S_{m,n}$ for all
$\al=k_{1}\al_{1}+k_{2}\al_{2}+k(1-\delta_{k_{1}k_{2},0})e_{2}$
with $k_{1},k_{2}\in\Z_{+}$ and $0\leq k\leq k_{1}+k_{2}$. By
Definition \ref{2.2}, the lemma holds.

{\bf Case 2.} $m>0$, $n<0$. Set
$$\al_{1}=(1-nm)e_{1}-n^{2}e_{2}, \
\al_{2}=m^{2}e_{1}+(1+nm)e_{2}.$$ Then $\al_{1},\al_{2}$ is a
$\Z$-basis of $\Z\times\Z$. Replacing $k$ by $-n$ and $m$ in
(\ref{ea2.1})-(\ref{ea2.3}) respectively, we have
$$\{L_{-nm+1,-n^{2}}, L_{-nm,-n^{2}}, L_{m^{2},mn+1},
L_{m^{2}-1,mn+1}\}\subseteq S_{m,n}.$$ Therefore,
$L_{\al}\subseteq {\mathcal S}_{m,n}$ for all
$\al=k_{1}\al_{1}+k_{2}\al_{2}+k(1-\delta_{k_{1}k_{2},0})e_{2}$
with $k_{1},k_{2}\in\Z_{+}$ and $0\leq k\leq k_{1}+k_{2}$.

For the case that $m<0$, $n>0$, set
$$\al_{1}=(1+nm)e_{1}+n^{2}e_{2}, \
\al_{2}=-m^{2}e_{1}+(1-nm)e_{2}.$$ Then
$$\{L_{nm+1,n^{2}}, L_{nm,n^{2}}, L_{-m^{2},-mn+1},
L_{-m^{2}-1,-mn+1}\}\subseteq {\mathcal S}_{m,n}.$$ For $m<0$,
$n<0$, set
$$\al_{1}=(-1-nm)e_{1}-n^{2}e_{2}, \
\al_{2}=-m^{2}e_{1}+(1-nm)e_{2}.$$ Then
$$\{L_{-nm-1,-n^{2}}, L_{-nm,-n^{2}}, L_{-m^{2},-mn+1},
L_{-m^{2}+1,-mn+1}\}\subseteq {\mathcal S}_{m,n}.$$ We can
similarly deduce the lemma for these two cases. \qed

\begin{lem}\label{2.4} Let $V=\bigoplus_{r,s\in\Z}V_{r,s}$ be an irreducible
${\mathcal L}$-module and $m,n\in\Z$ be such that $|m|\leq 1$ and
$|n|\geq 3$. Denote by ${\mathcal S}'_{m,n}$ the Lie subalgebra of
${\mathcal L}$ generated by $\{L_{m+i,n+j}|i=-1,0,1, j=-2,1,4\}$.
If there exists $0\neq v\in V$ such that $S'_{m,n}\cdot v=0$, then
$V$ is a $GHW$ ${\mathcal L}$-module.
\end{lem}

\pf We only consider the case that $m=1$, $n\geq 3$. For other
cases, the proof is similar. Set
$$\al_{1}=(n-1)e_{1}+n^{2}e_{2}, \
\al_{2}=e_{1}+(n+1)e_{2}.$$ Then $\al_{1},\al_{2}$ is a $\Z$-basis
of $\Z\times\Z$. Similar to the proof of Lemma \ref{2.3}, we
have$$ \{L_{n-1,n^{2}}, L_{n,n^{2}}, L_{2,n+1},
L_{1,n+1}\}\subseteq {\mathcal S}'_{1,n}.$$ So $L_{\al}\subseteq
{\mathcal S}'_{m,n}$  for all
$\al=k_{1}\al_{1}+k_{2}\al_{2}+k(1-\delta_{k_{1}k_{2},0})e_{2}$
with $k_{1},k_{2}\in\Z_{+}$ and $0\leq k\leq k_{1}+k_{2}$.
Therefore, the lemma holds.

\begin{definition}\label{2.5}
Let $V=\bigoplus_{r,s\in\Z}V_{r,s}$ be a ${\mathcal L}$-module. If
there exists a positive integer $N$ such that $\dim V_{r,s}\leq N$
for all $r,s\in\Z$, then $V$ is called uniformly bounded.
\end{definition}

We now have the first main result of the paper.

\begin{theorem}\label{2.6} Let $V=\bigoplus_{r,s\in\Z}V_{r,s}$ be an irreducible
${\mathcal L}$-module such that $V$ is not uniformly bounded, then
$V$ is a $GHW$ module.
\end{theorem}

\pf Suppose that $V$ is not a $GHW$ module. Let $m,n\in\Z$.

{\bf Case 1.} $|m|\geq 2$, $n\neq 0$. Let
$S_{-m,-n}=\{L_{-m+i,-n+j}|i=1,0,-1, j=3,0,-3\}$ be the same as in
Lemma \ref{2.3}. Then by Lemma \ref{2.3}, for any $0\neq v\in
V_{m,n}$, we have $S_{-m,-n}\cdot v\neq\{0\}$. By the assumption
\ref{ea22.1},
$$(\sum\limits_{i=-1,0,1}(\sum\limits_{j=-3,0,3}L_{-m+i,-n+j}))|_{V_{m,n}}:
V_{m,n}\longrightarrow
\bigoplus_{i=-1,0,1}(\bigoplus\limits_{j=-3}^{5}V_{i,j})$$ is an
injection.
Therefore,
\begin{equation}\label{ea2.7}\dim V_{m,n}\leq
\sum\limits_{i=-1,0,1}(\sum\limits_{j=-3}^{5}\dim V_{i,j}).
\end{equation}

 {\bf Case 2.} $|m|=1, n\neq 0$. Then by Lemma \ref{2.4}, for any $0\neq v\in V_{m,n}$, $S'_{-m,-n}\cdot
v\neq\{0\}$. Similar to the proof for case 1, we have
\begin{equation}\label{ea2.8}\dim V_{m,n}\leq
\sum\limits_{i=-1,0,1}(\sum\limits_{j=-2}^{6}\dim
V_{i,j}).\end{equation}

{\bf Case 3.} $m=0, n\neq 0$ or $m\neq 0, n=0$.  If $m=0,n\neq 0$,
by Lemma \ref{2.4},  for any $0\neq v\in V_{0,n}$, we have
$S'_{1,-n}\cdot v\neq\{0\}$. Therefore,
\begin{equation}\label{ea2.9}\dim V_{0,n}\leq
\sum\limits_{i=0,1,2}(\sum\limits_{j=-2}^{6}\dim
V_{i,j}).\end{equation}
 If $m\neq 0$, $n=0$, then $n-1=-1$ and
similarly we have
\begin{equation}\label{ea2.10}\dim V_{m,0}\leq
\sum\limits_{i=-1,0,1}(\sum\limits_{j=-4}^{4}\dim V_{i,j}), \ \
{\rm if} \ \ |m|\geq 2,
\end{equation}
\begin{equation}\label{ea2.11}\dim V_{m,0}\leq
\sum\limits_{i=-1,0,1}(\sum\limits_{j=-3}^{5}\dim V_{i,j}), \ \
{\rm if} \ \ |m|=1.
\end{equation}
It follows from
(\ref{ea2.7})-(\ref{ea2.11}) that $V$ is uniformly bounded, a
contradiction. \qed


\section{\bf Uniformly Bounded Modules}
\label{Preliminaries}
\label{sub3-c-related}

In this section, we discuss the  ${\mathcal L}$-modules
$V=\bigoplus_{r,s\in\Z}V_{r,s}$ such that $\dim V_{r,s}<\infty$
and
\begin{equation}\label{ea3.1}L_{a,b}V_{r,s}\subseteq
V_{r+a,s+b+1}+V_{r+a,s+b}.
\end{equation}

\begin{lem}\label{3.1} Let $V=\bigoplus_{r,s\in\Z}V_{r,s}$ be a
 ${\mathcal L}$-module satisfying (\ref{ea3.1}). Then $c$ acts
trivially on $V$.
\end{lem}

\pf By the definition of ${\mathcal L}$, we have the following Lie
bracket relation:
\begin{eqnarray*}
& [L_{\alpha_{1},\beta_{1}},L_{\alpha_{2},\beta_{2}}]&
=(\alpha_{2}-\alpha_{1})L_{\alpha_{1}+\alpha_{2},\beta_{1}+\beta_{2}+1}
+(\beta_{2}-\beta_{1})L_{\alpha_{1}+\alpha_{2},\beta_{1}+\beta_{2}}
\\
& &+\frac{1}{12}\delta_{\alpha_{1}+\alpha_{2},
0}[\delta_{\beta_{1}+\beta_{2}, -3}\alpha_{1}^{3}
+3\delta_{\beta_{1}+\beta_{2}, -2}(\beta_{1}+1)\alpha_{1}^{2}
\\
& &+3\delta_{\beta_{1}+\beta_{2},
-1}\beta_{1}(\beta_{1}+1)\alpha_{1}+\delta_{\beta_{1}+\beta_{2},
0}\beta_{1}(\beta_{1}^{2}-1)]c.
\end{eqnarray*}
Let $\al_{1}=-\al_{2}$, $\beta_{1}+\beta_{2}=k$, $k=-3,-2,-1,0$
respectively, then by (\ref{ea3.1}) and the above relation, we
have
$$c\cdot V_{m,n}\subseteq \bigcap_{k=-3}^{0}(\bigoplus_{k\leq i\leq
k+2}V_{m,n+i})=\{0\},$$ for all $m,n\in\Z$. The lemma is proved.

\qed

For $m,n\in \Z$ and $L_{m,n}\in{\mathcal L}$, by the assumption
(\ref{ea3.1}), for each $v\in V_{r,s}$,
$$L_{m,n}\cdot v=v_{1}+v_{2},$$
where $v_{1}\in V_{r+m,s+n+1}$, $v_{2}\in V_{r+m,s+n}$,
$r,s\in\Z$. Define $P_{m,n}: V_{r,s}\rightarrow V_{r+m,s+n+1}$ and
$Q_{m,n}: V_{r,s}\rightarrow V_{r+m,s+n}$ by
$$P_{m,n}\cdot v=v_{1}, \ \ Q_{m,n}\cdot v=v_{2}.$$ Then
$$L_{m,n}=P_{m,n}+Q_{m,n}.$$

\begin{lem}\label{l3.2} Let $V=\bigoplus_{r,s\in\Z}V_{r,s}$ be a
 ${\mathcal L}$-module satisfying (\ref{ea3.1}). For $h,k,r,s\in\Z$,  we have
\begin{equation}\label{eea3.1}
[P_{h,k}, P_{r,s}]=(r-h)P_{h+r,k+s+1},\end{equation}
\begin{equation}\label{eea3.2}
[P_{h,k},Q_{r,s}]+[Q_{h,k},P_{r,s}]=(r-h)Q_{h+r,k+s+1}+(s-k)P_{h+r,k+s},\end{equation}
\begin{equation}\label{eea3.3}[Q_{h,k},Q_{r,s}]=(s-k)Q_{h+r,k+s}.
\end{equation}
In particular, $${\mathcal P}_{-1}=span\{ \ P_{m,-1} \ | \ m\in
\Z\}$$ and
$${\mathcal Q}_{0}=span\{ \ Q_{0,n} \ | \ n\in \Z\}$$
are  centerless Virasoro Lie algebras and for each $m,n\in\Z$,
$$
V^{L}_{m}=\bigoplus_{s\in\Z}V_{m,s}, \ \
V^{R}_{n}=\bigoplus_{r\in\Z}V_{r,n}$$ are ${\mathcal Q}_{0}$ and
${\mathcal P}_{-1}$ modules respectively.
\end{lem}

In the following sections of the paper, we discuss the
indecomposable ${\mathcal L}$-modules
$V=\bigoplus_{r,s\in\Z}V_{r,s}$ satisfying (\ref{ea3.1}) and $\dim
V_{m,n}\leq 1$, for all $m,n\in\Z$. By Lemma \ref{l3.2}, we can
assume that
\begin{equation}\label{ea3.4}P_{0,-1}u=(\lambda+m)u, \ \
Q_{0,0}v=(\mu+n)v, \end{equation} for $u\in V_{m,-1}$, $v\in
V_{0,n}$, where $\lambda,\mu\in\C$ are two fixed complex numbers.

For convenience, denote
$$
P_{h,k}v_{m,n}=f_{h,k}(m,n)v_{h+m,n+k+1}, \ \
Q_{h,k}v_{m,n}=g_{h,k}(m,n)v_{m+h,n+k},$$ for $m,n,h,k\in\Z$. Then
\begin{equation}\label{e3.4}
\begin{array}{ll}
& f_{h,k}(m+r,n+s+1)f_{r,s}(m,n)-f_{r,s}(h+m,k+n+1)f_{h,k}(m,n)
\\
&=(r-h)f_{r+h,s+k+1}(m, n)
\end{array}\end{equation}
\\
\begin{equation}\label{e3.5}
\begin{array}{ll}
&f_{h,k}(r+m,s+n)g_{r,s}(m,n)+g_{h,k}(m+r,s+n+1)f_{r,s}
(m, n) \\
&-f_{r,s}(h+m,k+n)g_{h,k}(m,n)-g_{r,s}(m+h,k+n+1)f_{h,k}(m, n)
\\
&=(r-h)g_{h+r,k+s+1}(m,n)+(s-k)f_{h+r,k+s}( m, n),
\end{array}\end{equation}
\\
\begin{equation}\label{e3.6}
\begin{array}{ll}
&  g_{h,k}(
m+r,s+n)g_{r,s}(m,n)-g_{r,s}(m+h, n+k)g_{h,k}(m,n)\\
&=(s-k)g_{h+r, k+s}(m, n).
\end{array}\end{equation}

By (\ref{ea3.4}), we have
\begin{equation}\label{eea3.4}f_{0,-1}(m,-1)=\lambda+m,\ \ g_{0,0}(0,n)=\mu+n.
\end{equation}

\begin{lem}\label{l3.3} We have
\begin{equation}\label{3.10}f_{0,-1}(m,n)=\lambda+m\end{equation}
\begin{equation}\label{eaa3.11}
g_{0,0}(m,n)=\mu+n, \end{equation} for all $m,n\in\Z$.
\end{lem}

\pf Let $m=r=s=k=0$ in (\ref{e3.6}), then
\begin{equation}\label{ea3.11}
g_{h,0}(0,n)g_{0,0}(0,n)=g_{0,0}(h,n)g_{h,0}(0,n), \end{equation}
for all $n,h\in\Z$. Let $k=-1, $ $s=r=h=0$ in (\ref{e3.4}), then
$$f_{0,-1}(m,n+1)f_{0,0}(m,n)=f_{0,0}(m,n)f_{0,-1}(m,n).$$

If for each $h\in\Z$, there exists $n\in\Z$ such that
$g_{h,0}(0,n)\neq 0$, then (\ref{eaa3.11}) follows from
(\ref{ea3.11}).

Suppose there exists $h\in\Z$ such that $g_{h,0}(0,n)=0$ for all
$n\in\Z$. Let $m=r=k=0$ and $s\neq 0$ in (\ref{e3.6}), we have
$$
g_{h,s}(0,n)=0,$$ for all $s,n\in\Z$.  Let $r=-h$ and $s=k=m=0$ in
(\ref{e3.6}), then
$$g_{h,0}(-h,n)g_{-h,0}(0,n)=0,$$
for all $n\in\Z$.
 Let $r=-h$ and
$k=m=0$ in (\ref{e3.6}),
 we have
$$g_{h,0}(-h,s+n)g_{-h,s}(0,n)=sg_{0,s}(0,n).$$
So for $s\neq 0$ and $n\in\Z$, if $g_{h,0}(-h,s+n)=0$, then
$g_{0,s}(0,n)=0$. If $g_{h,0}(-h,s+n)\neq 0$, then
$g_{-h,0}(0,s+n)=0$. Since $g_{h,-s}(0,s+n)=0$, let $m=s=0, r=-h,
k=-s$ and replace $n$ by $n+s$ in (\ref{e3.6}), then we have
$$g_{0,-s}(0,s+n)=0.$$
Therefore for $s,n\in\Z$, $s\neq 0$, we have
$$g_{0,s}(0,n)g_{0,-s}(0,s+n)=0.$$
Let $m=h=r=0$, $k=-s$ in (\ref{e3.6}), then
$$g_{0,0}(0,n)=0$$
for all $n\in\Z$, which is not true. \qed
\\

We now have the second main result of the paper.

\begin{theorem}\label{t3.1} Let $V=\bigoplus_{r,s\in\Z}{\C}v_{r,s}$ be
an indecomposable ${\mathcal L}$-module satisfying (\ref{ea3.1})
and (\ref{ea3.4}). Then $V$ satisfies one of the following
situations:

(1) $A_{a,\lambda,\mu}:$
$L_{r,s}v_{m,n}=(ar+\lambda+m)v_{m+r,n+s+1}+(as+\mu+n)v_{m+r,n+s}$,
where $\lambda, \mu, a\in\mathbb{C}$;

(2) $A_{0,\lambda,\mu}:$
$L_{r,s}v_{m,n}=\displaystyle\frac{\mu+n}{\mu+n+s+1}(\lambda+m)v_{m+r,n+s+1}
+\displaystyle\frac{\lambda+m}{\lambda+m+r}(\mu+n)v_{m+r,n+s},$
where $\lambda, \mu\not\in\mathbb{Z}$;

(3) $A_{1,\lambda,\mu}:$
$L_{r,s}v_{m,n}=\displaystyle\frac{\mu+n+s+1}{\mu+n}(r+\lambda+m)v_{m+r,n+s+1}
+\displaystyle\frac{\lambda+m+r}{\lambda+m}(s+\mu+n)v_{m+r,n+s}$,
where $\lambda, \mu\not\in\mathbb{Z}$;

(4) $A_{1,0,\lambda,\mu}:$
$L_{r,s}v_{m,n}=(r+\lambda+m)v_{m+r,n+s+1}
+\displaystyle\frac{\lambda+m+r}{\lambda+m}(\mu+n)v_{m+r,n+s}$,
where $\lambda \not\in\mathbb{Z}$;

(5) $B_{1,0,\lambda,\mu}:$
$L_{r,s}v_{m,n}=\displaystyle\frac{\mu+n}{\mu+n+s+1}(r+\lambda+m)v_{m+r,n+s+1}
+(\mu+n)v_{m+r,n+s}$, where $\mu\not\in\mathbb{Z}$;

(6) $A_{0,1,\lambda,\mu}:$
$L_{r,s}v_{m,n}=(\lambda+m)v_{m+r,n+s+1}
+\displaystyle\frac{\lambda+m}{\lambda+m+r}(s+\mu+n)v_{m+r,n+s}$,
where $\lambda\not\in\mathbb{Z}$;

(7) $B_{0,1,\lambda,\mu}:$
$L_{r,s}v_{m,n}=\displaystyle\frac{\mu+n+s+1}{\mu+n}(\lambda+m)v_{m+r,n+s+1}
+(s+\mu+n)v_{m+r,n+s}$, where $\mu\not\in\mathbb{Z}$.
\end{theorem}

We will prove the theorem through lemmas in the next two sections.
It is easy to check that (1)-(7) are indecomposable representations
of ${\mathcal L}$.

\section{\bf Proof of Theorem \ref{t3.1} for some cases}
\label{Preliminaries}
\label{sub4-c-related}

In this section  we  suppose that
\begin{equation}\label{5.1}
Q_{0,s}v_{m,n}=(a_{m}s+\mu+n)v_{m,n+s},
\end{equation}
\begin{equation}\label{5.2}
P_{r,-1}v_{m,n}=(b_{n}r+\lambda+m)v_{m+r,n},
\end{equation}
where $a_{m}, b_{n}\in\mathbb{C}$. In fact, if $\mu\neq0$ or
$\lambda\neq 0$, then $Q_{0,s}$ or $P_{r,-1}$ has to be the form of
(\ref{5.1}) or (\ref{5.2}), since
$V^{L}_{m}=\bigoplus_{n\in\Z}\mathbb{C}v_{m,n}$ and
$V^{R}_{n}=\bigoplus_{m\in\Z}\mathbb{C}v_{m,n}$ are ${\mathcal
Q}_{0}$ and ${\mathcal P}_{-1}$ modules respectively.

\begin{lem}\label{l3.4}
$a_{m}=a$ for all $m\in\mathbb{Z}$. That is,
$$g_{0,s}(m,n)=as+\mu+n,$$
for all $m\in\mathbb{Z}$.
\end{lem}

\pf By the fact that $[L_{0, -1}, L_{r,-1}]=rL_{r, -1}$ and
(\ref{eea3.4}), we have
\begin{equation}\label{5.0}
g_{0,-1}(m+r,n)f_{r,-1}(m,n)=g_{0,-1}(m,n)f_{r,-1}(m,n-1).
\end{equation}
So
\begin{equation}\label{5.3}
(\mu+n-a_{m+r})(b_{n}r+\lambda+m)=(\mu+n-a_{m})(b_{n-1}r+\lambda+m),
\end{equation}
for all $m,n,r\in\mathbb{Z}.$ Replacing $r$ by $-r$ and $m$ by $m+r$
in (\ref{5.3}), we have
\begin{equation}\label{5.4}
(\mu+n-a_{m})(-b_{n}r+\lambda+m+r)=(\mu+n-a_{m+r})(-b_{n-1}r+\lambda+m+r).
\end{equation}
By (\ref{5.3})and (\ref{5.4}) we have
$(b_{n}+b_{n-1}-1)r(a_{m+r}-a_{m})=0$. Then
$$b_{n}+b_{n-1}=1\; or \;a_{m+r}=a_{m},$$
for all $n,m,r\in\Z$. If $b_{n}+b_{n-1}=1$ for  all $n\in\Z$, then
$$b_{2k}=b_{0},\quad  b_{2k+1}=1-b_{0},\quad\quad\;k\in\Z.$$
Let $n=2k$ and  $n=2k+1$ in (\ref{5.3}) respectively,  then we
have
\begin{equation}\label{1}
(2b_{0}-1)r(\mu+2k-a_{m})=(b_{0}r+\lambda+m)(a_{m+r}-a_{m}),
\end{equation}
\begin{equation}\label{2}
(2b_{0}-1)r(\mu+2k+1-a_{m+r})=(b_{0}r+\lambda+m)(a_{m}-a_{m+r}),
\end{equation}
for all $k, m, r\in\Z$. By (\ref{1})and (\ref{2}),  we have
$$(2b_{0}-1)(2\mu+4k-1-a_{m}-a_{m+r})=0, $$ for all $k\in\Z$.
There always exists $k\in\Z$ such that
$(2\mu+4k-1-a_{m}-a_{m+r})\neq 0$ for all $m,r\in\Z$. Therefore,
$b_{0}=\dis\frac{1}{2}$. Then $b_{n}=\dis\frac{1}{2}$ for all
$n\in\Z$. By (\ref{5.3}), we have
\begin{equation}\label{3}
(\frac{1}{2}r+\lambda+m)(a_{m}-a_{m+r})=0,
\end{equation}
for all $m,r\in\Z$. Let $m=0$ in (\ref{3}), then
$$(\frac{1}{2}r+\lambda)(a_{0}-a_{r})=0.$$
If $\lambda\not\in\dis\frac{1}{2}\Z$, then $a_{r}=a_{0}$ for all
$r\in\Z$. If $\lambda\in\dis\frac{1}{2}\Z$, by (\ref{3}), we have
$$a_{r}=a_{0},$$
for all $r\neq-2\lambda$. Let $m=-2\lambda,  r\neq 0$ and  $r\neq
2\lambda$ in (\ref{3}), then we have $a_{-2\lambda}=a_{0}$. So
$a_{m}=a$ for all $m\in\Z$. \qed

 From the relation that
$[L_{r,-1}, L_{0,s}]=-rL_{r,s}+(s+1)L_{r,s-1}$, we have
\begin{equation}\label{5.15-}
(b_{n+s+1}r+\lambda+m)f_{0,s}(m,n)-(b_{n}r+\lambda+m)f_{0,s}(m+r,n)=-rf_{r,s}(m,n),
\end{equation}
\begin{equation}\label{5.16-}
\begin{array}{ll}
&  (b_{n+s+1}-b_{n})r(as+\mu+n)+f_{0,s}(m,n)g_{r,-1}(m,n+s+1)\\
&-g_{r,-1}(m,n)f_{0,s}(m+r,n-1)=-rg_{r,s}(m,n)+(s+1)f_{r,s-1}(m,n),
\end{array}\end{equation}
\begin{equation}\label{5.17}
g_{r,-1}(m,n+s)(as+\mu+n)-g_{r,-1}(m,n)(as+\mu+n-1)=(s+1)g_{r,s-1}(m,n).
\end{equation}
By the fact that $[L_{0,-1}, L_{r,-1}]=rL_{r,-1}$, we have
\begin{equation}\label{5.9}
g_{r,-1}(m,n)(\mu+n-1-a)=g_{r,-1}(m,n-1)(\mu+n-a).
\end{equation}

By (\ref{5.3}), if $\mu+n-a\neq 0$ for all $n\in\mathbb{Z}$, then
$b_{n}=b$ for all $n\in\mathbb{Z}$.

\begin{lem}\label{l3.5}
Suppose that $\mu+n-a\neq0$ for all $n\in\mathbb{Z}$.

(1)  If $b\neq0,1$, then a=b and
\begin{equation}
\begin{cases}\nonumber
 & f_{r,s}(m,n)=ar+\lambda+m,\\
 & g_{r,s}(m,n)=as+\mu+n;
\end{cases}
\end{equation}


(2) If $b\in\{0,1\}$, then $a\in\{0,1\}$.
\end{lem}

\pf By  (\ref{5.15-})-(\ref{5.16-}), we have
\begin{equation}\label{5.15}
(br+\lambda+m)(f_{0,s}(m,n)-f_{0,s}(m+r,n))=-rf_{r,s}(m,n),
\end{equation}
\begin{equation}\label{5.16}
\begin{array}{ll}
&  f_{0,s}(m,n)g_{r,-1}(m,n+s+1)-g_{r,-1}(m,n)f_{0,s}(m+r,n-1)\\
&=-rg_{r,s}(m,n)+(s+1)f_{r,s-1}(m,n).
\end{array}\end{equation}
By (\ref{5.9}) and  $\mu+n-a\neq0$ for all $n\in\mathbb{Z}$, we
can assume that
$$\frac{g_{r,-1}(m,n)}{\mu+n-a}=\frac{g_{r,-1}(m,n-1)}{\mu+n-1-a}=k_{r,m},$$
for all $m,r,n\in\mathbb{Z}$. Note that $k_{0,m}=1$ for all
$m\in\mathbb{Z}$.
By (\ref{5.17}), we have
$$g_{r,s}(m,n)=k_{r,m}(as+\mu+n),\quad\quad\quad s\neq-2.$$
By (\ref{e3.6}), we have
\begin{equation}\label{5.10}
k_{r,m}k_{-r,m+r}=1,
\end{equation}
\begin{equation}\label{5.11}
k_{h,m+r}k_{r,m}=k_{r,m+h}k_{h,m}=k_{h+r,m},
\end{equation}
for all $r,m,h\in\mathbb{Z}$.  From (\ref{5.10}) we can see that
$k_{r,m}\neq 0$ for all $r,m\in\mathbb{Z}$. Let $h=0, k=-2, s=0$
in (\ref{e3.6}), we have $g_{r,-2}(m,n)=k_{r,m}(-2a+\mu+n).$  Then
\begin{equation}\label{5.12-}
g_{r,s}(m,n)=k_{r,m}(as+\mu+n),
\end{equation}
for all $r,s,m,n\in\Z$. In (\ref{e3.5}), let $h=0, r=0, k=-1$, then
we have
\begin{equation}\label{5.12}
f_{0,s}(m,n)(\mu+n+s+1-a)-(\mu+n-a)f_{0,s}(m,n-1)=(s+1)f_{0,s-1}(m,n).
\end{equation}
So
\begin{equation}\label{5.13}
f_{0,0}(m,n)(\mu+n+1-a)-f_{0,0}(m,n-1)(\mu+n-a)=\lambda+m,
\end{equation}
and
\begin{equation}\label{5.14}
(\mu+n+1-a)f_{0,0}(m,n)=n(\lambda+m)+(\mu+1-a)f_{0,0}(m,0).
\end{equation}
Letting $n=1,s=0$ in (\ref{5.16}),  we obtain
$$f_{0,0}(m,1)k_{r,m}(\mu+2-a)-k_{r,m}(\mu+1-a)f_{0,0}(m+r,0)=-rk_{r,m}(\mu+1)+br+\lambda+m.$$
By (\ref{5.14}), we have
\begin{equation}\label{5.18}
[\lambda+m+f_{0,0}(m,0)(\mu+1-a)]k_{r,m}-f_{0,0}(m+r,0)(\mu+1-a)k_{r,m}=br+\lambda+m-rk_{r,m}(\mu+1).
\end{equation}
So
\begin{equation}\label{5.19}
k_{r,m}(\mu+1-a)(f_{0,0}(m,0)-f_{0,0}(m+r,0))+k_{r,m}(\lambda+m+r\mu+r)=br+\lambda+m.
\end{equation}
Replace $r$ by $-r$ and $m$ by $m+r$ in (\ref{5.18}) respectively,
then
$$k_{-r,m+r}(\mu+1-a)(f_{0,0}(m+r,0)-f_{0,0}(m,0))+k_{-r,m+r}(\lambda+m-r\mu)=-br+\lambda+m+r.$$
Then by (\ref{5.10}), we have
$$k_{r,m}^{2}(-br+\lambda+m+r)-(2\lambda+2m+r)k_{r,m}+br+\lambda+m=0,$$
i.e.,
\begin{equation}\label{5.20}
(k_{r,m}-1)[(-br+\lambda+m+r)k_{r,m}-(br+\lambda+m)]=0.
\end{equation}

If $b=\dis\frac{1}{2}$, then $k_{r,m}=1$ for all $r,m\in\mathbb{Z}$.

Suppose $b\neq\dis\frac{1}{2}$ and there exist $r,m\in\mathbb{Z}$
such that $-br+\lambda+m+r\neq0$ and
$k_{r,m}=\dis\frac{br+\lambda+m}{-br+\lambda+m+r}\neq1$.  For
$h\neq-r, h\neq0$, if $k_{h,m+r}=1$, then by (\ref{5.11}) we have
$$k_{h+r,m}=k_{r,m}\neq1.$$
Then
$$k_{h+r,m}=\frac{b(h+r)+\lambda+m}{-b(h+r)+\lambda+m+r+h}=\frac{br+\lambda+m}{-br+\lambda+m+r}.$$
This implies that $b=\dis\frac{1}{2}$, a contradiction. So
$k_{h,m+r}\neq1$ for $h\neq0, h\neq-r$. Then by (\ref{5.11}),
$$\frac{bh+\lambda+m+r}{-bh+\lambda+m+r+h}\cdot\frac{br+\lambda+m}{-br+\lambda+m+r}
=\frac{b(h+r)+\lambda+m}{-b(h+r)+\lambda+m+r+h}.$$ This forces that
$b=0$ or $1$. By the assumption that $b\neq0,1$,  we have
$k_{r,m}=1$.
Then by (\ref{5.12-}), we have
\begin{equation}\label{5.22}
g_{r,s}(m,n)=as+\mu+n.
\end{equation}
 By (\ref{5.18}), we have
$$f_{0,0}(m+r,0)-f_{0,0}(m,0)=\frac{r(\mu+1-b)}{\mu+1-a}.$$
Then by (\ref{5.15}),
\begin{equation}\label{5.23}
f_{r,0}(m,0)=\frac{\mu+1-b}{\mu+1-a}(br+\lambda+m),\quad\quad\quad
r\neq0.
\end{equation}
Setting $h=r, r=s=k=n=0$ in (\ref{e3.5}),  we get
$$f_{r,0}(m,0)\mu+f_{0,0}(m,0)(\mu+1)-f_{0,0}(r+m,0)\mu-f_{r,0}(m,0)(\mu+1)=-r(a+\mu).$$
So
$$f_{r,0}(m,0)=f_{0,0}(m,0)-\frac{r\mu(\mu+1-b)}{\mu+1-a}+r(a+\mu).$$
Then by (\ref{5.23}), we have
\begin{equation}\label{5.24}
f_{0,0}(m,0)=(\lambda+m)\frac{\mu+1-b}{\mu+1-a}+\frac{r(b-b^{2}-a+a^{2})}{\mu+1-a},\;\;
r\neq 0.
\end{equation}
Therefore, $b-b^{2}-a+a^{2}=0.$ Then  $$a=b\;\; {\rm or} \;\;
a+b=1.$$
 If $a=b$, then $f_{0,0}(m,0)=\lambda+m$. By (\ref{5.14}),
\begin{equation}\label{5.25}
f_{0,0}(m,n)=\lambda+m.
\end{equation}
In (\ref{e3.5}), let $h=k=0$, then we have
\begin{equation}\label{5.26}
\begin{array}{ll}
&f_{0,0}(r+m,s+n)g_{r,s}(m,n)+g_{0,0}(m+r,s+n+1)f_{r,s}(m, n)
 \\
&-f_{r,s}(m,n)g_{0,0}(m,n)-g_{r,s}(m,n+1)f_{0,0}(m, n)
\\
&=rg_{r,s+1}(m,n)+sf_{r,s}( m, n).
\end{array}\end{equation}
By (\ref{5.22}) and (\ref{5.26}), we have
$$f_{r,s}(m,n)=ar+\lambda+m.$$
So we can deduce that
\begin{equation}
\begin{cases}\nonumber
 & f_{r,s}(m,n)=ar+\lambda+m,\\
 & g_{r,s}(m,n)=as+\mu+n.
\end{cases}
\end{equation}
If $a+b=1$, then
\begin{equation}\label{5.27}
f_{0,0}(m,n)=\frac{\mu+n+a}{\mu+n+1-a}(\lambda+m).
\end{equation}
By induction on $n$, we can deduce that
$$f_{0,1}(m,n)=\frac{(\mu+n+a)(\mu+n+1+a)}{(\mu+n+1-a)(\mu+n+2-a)}(\lambda+m)$$
and
\begin{equation}\label{5.28}
f_{r,1}(m,n)=\frac{(\mu+n+a)(\mu+n+1+a)}{(\mu+n+1-a)(\mu+n+2-a)}((1-a)r+\lambda+m).
\end{equation}
On the other hand, by (\ref{5.25}) and (\ref{5.26}), we have
\begin{equation}\label{5.29}
f_{r,1}(m,n)=ra+\frac{(\mu+n+a)}{(\mu+n+2-a)}[r(1-2a)+\frac{(\mu+n+1+a)}{(\mu+n+1-a)}(\lambda+m)].
\end{equation}
Then by (\ref{5.28}) and (\ref{5.29}), we have $$a=0, 1\; {\rm or}\;
\frac{1}{2}.$$ Then $b=1,0$ or $\dis\frac{1}{2}$. By the assumption,
$a\neq 0, 1$, so $a=\dis\frac{1}{2}=b$. This also means that
\begin{equation}
\begin{cases}\nonumber
 & f_{r,s}(m,n)=ar+\lambda+m,\\
 & g_{r,s}(m,n)=as+\mu+n.
\end{cases}
\end{equation}
(2) follows from the proof of (1). \qed

\begin{lem}\label{l3.6}
Suppose that $\mu+n-a\neq0$ for all $n\in\mathbb{Z}$ and $a=b=1$.
Then we have
\begin{equation}\label{6.1}
\begin{cases}
 & f_{r,s}(m,n)=r+\lambda+m,\\
 & g_{r,s}(m,n)=s+\mu+n,
\end{cases}
\end{equation}
or
\begin{equation}\label{6.2}
\begin{cases}
 & f_{r,s}(m,n)=\displaystyle\frac{\mu+n+s+1}{\mu+n}(r+\lambda+m),\\
 &
 g_{r,s}(m,n)=\displaystyle\frac{\lambda+m+r}{\lambda+m}(s+\mu+n),
\end{cases}\quad\quad\quad\lambda\not\in\mathbb{Z}.
\end{equation}
\end{lem}

\pf {\bf Case 1.}  $\lambda\not\in\mathbb{Z}$. By the proof of (1)
of Lemma \ref{l3.5} ( see (\ref{5.20}) ), we know that
$$k_{r,m}=1\;\; or \;\; k_{r, m}=\frac{r+\lambda+m}{\lambda+m}.$$
If $k_{r,m}=1$, then by (\ref{5.12-}), $g_{r,s}(m,n)=s+\mu+n$.  By
(\ref{5.24}) and (\ref{5.14}), we have
$$f_{0,0}(m,n)=\lambda+m.$$
Then by (\ref{5.26}), we have
$$f_{r,s}(m,n)=r+\lambda+m.$$
So we have (\ref{6.1}). If $k_{r,
m}=\displaystyle\frac{r+\lambda+m}{\lambda+m}$, then
\begin{equation}\label{5.30}
g_{r,s}(m,n)=\displaystyle\frac{r+\lambda+m}{\lambda+m}(s+\mu+n).
\end{equation}
By (\ref{5.19}), we have
$$f_{0,0}(m+r,0)-f_{0,0}(m,0)=\frac{\mu+1}{\mu}r.$$
By (\ref{5.15}), we obtain
$$f_{r,0}(m,0)=\frac{\mu+1}{\mu}(r+\lambda+m),\quad\quad\quad r\neq0.$$
In (\ref{e3.5}),  let $n=0, k=r=s=0$, then we can deduce that
$$f_{0,0}(m,0)=\frac{\mu+1}{\mu}(\lambda+m).$$
Then by (\ref{5.14}), we have
\begin{equation}\label{5.31}
f_{0,0}(m,n)=\displaystyle\frac{\mu+n+1}{\mu+n}(\lambda+m).
\end{equation}
By (\ref{5.26}), (\ref{5.30}) and (\ref{5.31}), we have (\ref{6.2}).

{\bf Case 2.}  $\lambda\in\mathbb{Z}$. By (\ref{5.15}), we have
\begin{equation}\label{5.32}
f_{r,s}(-\lambda-r,n)=0,\quad\quad\quad r\neq0.
\end{equation}
In (\ref{e3.4}), let $h=-r, k=-1, m=-\lambda$, then by
(\ref{5.32}), we have
\begin{equation}\label{5.33}
f_{0,s}(-\lambda,n)=0.
\end{equation}
By (\ref{5.20}),  we have $$k_{r,-\lambda-r}=1,\quad\quad\quad
r\neq0.$$ Then by $k_{r,-\lambda}=1$ for all $r\in\mathbb{Z}$ and
(\ref{5.11}), we have
$$k_{h+r,-\lambda-r}=k_{h,-\lambda}k_{r,-\lambda-r}=1,\;\;\;
\forall\; h, r\in\mathbb{Z}.$$ So
 $$k_{r,m}=1,\;\;\;\forall\; r,m\in\mathbb{Z}.$$
Then similar to the proof above, we have (\ref{6.1}). \qed

\begin{lem}\label{l3.7}
Suppose that $\mu+n-a\neq0$ for all $n\in\mathbb{Z}$ and $b=1,
a=0$. Then
\begin{equation}\label{7.1}
\begin{cases}
 & f_{r,s}(m,n)=r+\lambda+m,\\
 & g_{r,s}(m,n)=\displaystyle\frac{\lambda+m+r}{\lambda+m}(\mu+n),
\end{cases}\quad\quad\quad\lambda\not\in\mathbb{Z},
\end{equation}
or
\begin{equation}\label{7.2}
\begin{cases}
 & f_{r,s}(m,n)=\displaystyle\frac{\mu+n}{\mu+n+s+1}(r+\lambda+m),\\
 & g_{r,s}(m,n)=\mu+n.
\end{cases}
\end{equation}
\end{lem}

\pf {\bf Case 1.} $\lambda\in\mathbb{Z}$. As the proof in (2) of
Lemma \ref{l3.6}, we have
$$k_{r,m}=1, $$
and
\begin{equation}\label{5.34}
g_{r,s}(m,n)=\mu+n.
\end{equation}
By (\ref{5.24})and  (\ref{5.14}), we have
$$f_{0,0}(m,n)=\frac{\mu+n}{\mu+n+1}(\lambda+m).$$
Then by (\ref{5.26}) and (\ref{5.34}), we get (\ref{7.2}).

{\bf Case 2.} $\lambda\not\in\mathbb{Z}$. Then
$$k_{r,m}=1 \; or \; k_{r,m}=\frac{\lambda+m+r}{\lambda+m}.$$
If $k_{r,m}=1$,  then we have (\ref{7.2}). If
$k_{r,m}=\displaystyle\frac{\lambda+m+r}{\lambda+m}$, then
\begin{equation}\label{5.36}
g_{r,s}(m,n)=\displaystyle\frac{\lambda+m+r}{\lambda+m}(\mu+n).
\end{equation}
Then we can deduce that
$$f_{0,0}(m,n)=\lambda+m.$$
By (\ref{5.36}) and (\ref{5.26}), we have (\ref{7.1}). \qed

\begin{lem}\label{l3.8}
Suppose that $\mu+n-a\not\in\mathbb{Z}$ for all $n\in\mathbb{Z}$
and $b=0, a=1$. Then
\begin{equation}\label{8.1}
\begin{cases}
 & f_{r,s}(m,n)=\displaystyle\frac{\mu+n+s+1}{\mu+n}(\lambda+m),\\
 & g_{r,s}(m,n)=s+\mu+n,
\end{cases}
\end{equation}
or
\begin{equation}\label{8.2}
\begin{cases}
 & f_{r,s}(m,n)=\lambda+m,\\
 & g_{r,s}(m,n)=\displaystyle\frac{\lambda+m}{\lambda+m+r}(s+\mu+n),
\end{cases}\quad\quad\quad\lambda\not\in\mathbb{Z}.
\end{equation}
\end{lem}

\pf {\bf Case 1.} $\lambda\in\mathbb{Z}$. As the proof above  we
have
$$k_{r,m}=1$$ and
\begin{equation}\label{5.38}
g_{r,s}(m,n)=s+\mu+n.
\end{equation}
By (\ref{5.19}), we have
$$f_{0,0}(m+r,0)-f_{0,0}(m,0)=\frac{\mu+1}{\mu}r.$$
Then by (\ref{5.15}), we obtain
$$f_{r,0}(m,0)=\frac{\mu+1}{\mu}(\lambda+m),\quad\quad\quad r\neq0.$$
Let $n=0, k=r=s=0$ in (\ref{e3.5}), then we have
$$f_{0,0}(m,0)=\frac{\mu+1}{\mu}(\lambda+m).$$
Then by (\ref{5.14}), we get
$$f_{0,0}(m,n)=\frac{\mu+n+1}{\mu+n}(\lambda+m).$$
Then by  (\ref{5.38}) and (\ref{5.26}), we obtain (\ref{8.1}).

{\bf Case 2.} If $\lambda\not\in\mathbb{Z}$, then
$$k_{r,m}=1\;\; {\rm or} \;\; k_{r, m}=\frac{\lambda+m}{\lambda+m+r}.$$
If $k_{r,m}=1$, then $g_{r,s}(m,n)=s+\mu+n.$ Using (\ref{5.24}) and
(\ref{5.14}), we have
$$f_{0,0}(m,n)=\frac{\mu+n+1}{\mu+n}(\lambda+m).$$
By  (\ref{5.26}), we have
$$f_{r,s}(m,n)=\frac{\mu+s+n+1}{\mu+n}(\lambda+m).$$
Then we get (\ref{8.1}).
 If $k_{r, m}=\displaystyle\frac{\lambda+m}{\lambda+m+r}$, then
$g_{r,s}(m,n)=\displaystyle\frac{\lambda+m}{\lambda+m+r}(s+\mu+n).$
By (\ref{5.19}), we have $$f_{0,0}(m+r, 0)-f_{0,0}(m, 0)=r.$$  By
(\ref{5.15}), we obtain $$f_{r,0}(m, 0)=\lambda+m,\;\;\;r\neq0.$$
Let $k=r=s=n=0, h\neq 0$ in (\ref{e3.5}), then we have
$$f_{0,0}(m, 0)=\lambda+m.$$
Then by (\ref{5.14}), we have $$f_{0,0}(m, n)=\lambda+m.$$ Using
(\ref{5.26}),  we deduce $$f_{r,s}(m, n)=\lambda+m.$$ So we have
(\ref{8.2}). \qed

Similarly we have the following lemma.

\begin{lem}\label{l3.9}
Suppose that $\mu+n-a\neq0$ for all $n\in\mathbb{Z}$ and $b=a=0$.
Then
\begin{equation}
\begin{cases}\nonumber
 & f_{r,s}(m,n)=\lambda+m,\\
 & g_{r,s}(m,n)=\mu+n,
\end{cases}
\end{equation}
or
\begin{equation}
\begin{cases}\nonumber
 & f_{r,s}(m,n)=\displaystyle\frac{\mu+n}{\mu+n+s+1}(\lambda+m),\\
 & g_{r,s}(m,n)=\displaystyle\frac{\lambda+m}{\lambda+m+r}(\mu+n),
\end{cases}\quad\quad\quad\lambda\not\in\mathbb{Z}.
\end{equation}
\end{lem}

\begin{lem}\label{l3.10}
Suppose that $\mu-a\in\mathbb{Z}$. Then we also have $b_{n}=b$ for
all $n\in\mathbb{Z}$. Furthermore, we have
\begin{equation}\label{10.1}
\begin{cases}
 & f_{r,s}(m,n)=ar+\lambda+m,\\
 & g_{r,s}(m,n)=as+\mu+n,
\end{cases}
\end{equation}
or
\begin{equation}\label{10.2}
\begin{cases}
 & f_{r,s}(m,n)=r+\lambda+m,\\
 & g_{r,s}(m,n)=\displaystyle\frac{\lambda+m+r}{\lambda+m}(\mu+n),
\end{cases}\quad\quad\quad\lambda\not\in\mathbb{Z},
\end{equation}
or
\begin{equation}\label{10.3}
\begin{cases}
 & f_{r,s}(m,n)=\lambda+m,\\
 & g_{r,s}(m,n)=\displaystyle\frac{\lambda+m}{\lambda+m+r}(s+\mu+n),
\end{cases}\quad\quad\quad\lambda\not\in\mathbb{Z}.
\end{equation}
\end{lem}

\pf Let $n_{0}\in\mathbb{Z}$ be such that $\mu+n_{0}-a=0.$  By
(\ref{5.0}), we have
$$(b_{n}r+\lambda+m)(\mu+n-a)=(b_{n-1}r+\lambda+m)(\mu+n-a).$$
So
\begin{equation}\nonumber
b_{n}=
\begin{cases}
 b_{n_{0}}, & n\geq n_{0},\\
 b_{n_{0}-1},& n\leq n_{0}-1.
\end{cases}
\end{equation}
Then by (\ref{5.14}), we have
$$f_{0,0}(m,n)=\lambda+m, \quad n\neq n_{0}-1.$$ Assume
$f_{0,0}(m,n_{0}-1)=\lambda+m+f(m)$ for all $m\in\Z$. Then
\begin{equation}\label{5.021}
f_{0,0}(m,n)=\lambda+m+\delta_{n, n_{0}-1}f(m).
\end{equation}
Since $\mu+n-a\neq0$ for all $n\neq n_{0}$, by (\ref{5.9}) we can
assume that
$$g_{r,-1}(m,n)=k_{r,m}(\mu+n-a),\quad for\;\;\; n\geq n_{0},$$
$$g_{r,-1}(m,n)=k'_{r,m}(\mu+n-a),\quad for\; n\leq n_{0}-1.$$
Similar to  the proof of Lemma \ref{l3.5},  we have
\begin{equation}\label{5.016}
k_{r,m}k_{-r,m+r}=1, \;\;\;
k_{h,m+r}k_{r,m}=k_{r,m+h}k_{h,m}=k_{h+r,m},
\end{equation}
\begin{equation}\label{5.017}
k'_{r,m}k'_{-r,m+r}=1,\;\;\;
k'_{h,m+r}k'_{r,m}=k'_{r,m+h}k'_{h,m}=k'_{h+r,m},
\end{equation}
for all $r,m,h\in\mathbb{Z}$. By (\ref{5.17}), we have
$$(s+2)g_{r,s}(m,n)=(as+\mu+n+a)g_{r,-1}(m,s+n+1)-(as+\mu+n+a-1)g_{r,-1}(m,n).$$
Then by (\ref{e3.6}), we have
\begin{eqnarray*}
&&g_{r,s}(m,n)=k_{r,m}(as+\mu+n),\quad  n\geq n_{0}, s\geq-1;
\\
&&g_{r,s}(m,n)=\frac{1}{s+2}[(as+\mu+n+a)(\mu+s+n+1-a)k'_{r,m}
\\
&&\quad\quad\quad\quad\quad -(as+\mu+n+a-1)(\mu+n-a)k_{r,m}],\quad
 n\geq n_{0}, s\leq -3;
\\
&&g_{r,s}(m,n)=\frac{1}{s+2}[(as+\mu+n+a)(\mu+s+n+1-a)k_{r,m}
\\
&&\quad\quad\quad\quad\quad
-(as+\mu+n+a-1)(\mu+n-a)k'_{r,m}],\quad  n\leq n_{0}-1, s\geq 0;
\\
&&g_{r,s}(m,n)=k'_{r,m}(as+\mu+n),\quad  n\leq n_{0}-1, s\leq -3.
\end{eqnarray*}
Let $h=0, k=-3, s=1$ in (\ref{e3.6}), then
$$g_{r,-2}(m,n)=\frac{1}{4}[(-3a+\mu+n+1)g_{r,1}(m,n)-(-3a+\mu+n)g_{r,1}(m,n-3)].$$
So we can deduce that
\begin{eqnarray*}
&&g_{r,-2}(m,n)=k_{r,m}(-2a+\mu+n),\quad n\geq n_{0}+3;
\\
&&g_{r,-2}(m,n)=\frac{1}{12}[3(-3a+\mu+n)(2\mu+2n+a-1)+(2a+\mu+n)(\mu+n+2-a)]k_{r,m}
\\
&&-\frac{1}{12}[3(-3a+\mu+n)(2\mu+2n+a-4)+(2a+\mu+n-1)(\mu+n-a)]k'_{r,m},\;
n\leq n_{0}-1;
\\
&&g_{r,-2}(m,n_{0})=\frac{1}{2}ak_{r,m};
\\
&&g_{r,-2}(m,n_{0}+1)=\frac{1}{2}[(1-a)(2a+1)k_{r,m}+(a-1)(2a-1)k'_{r,m}];
\\
&&g_{r,-2}(m,n_{0}+2)=\frac{1}{6}[(a+2)(-3a+5)k_{r,m}+(a-1)(3a-2)k'_{r,m}].
\end{eqnarray*}
If $a=1$, letting  $k=0, s=-2, n=n_{0}$ in (\ref{e3.6}),  we have
$$k'_{h,m+r}k_{r,m}+k_{r,m+h}k_{h,m}=2k_{r+h,m}.$$
Let $r=0$ and note $k'_{0,m}=1=k_{0,m}$ for all $m\in\Z$, then we
have
$$k'_{h,m}=k_{h,m}.$$
 If $a=0$, let  $k=-3,
h=0, s=-2, n=n_{0}+2$ in (\ref{e3.6}), then we obtain
$$g_{0,-3}(m+r,n_{0})g_{r,-2}(m,n_{0}+2)-g_{r,-2}(m,n_{0}-1)g_{0,-3}(m+r,n_{0}+2)=g_{r,-5}(m,n_{0}+2).$$
Then we  can also get  $k'_{r,m}=k_{r,m}$ for all $r,m\in\Z$.
Therefore, if $a=0$ or 1, we have $$k'_{r,m}=k_{r,m}, \;\;\forall\;
r,m\in\Z.$$
 Letting  $n=n_{0}+1, s=0$ in
(\ref{5.16-}),  we have
\begin{equation}\label{5.018}
(\lambda+m+ar)k_{r,m}=b_{n_{0}}r+\lambda+m.
\end{equation}
Letting $n=n_{0}-3, s=0$ in (\ref{5.16-}), we obtain
\begin{equation}\label{5.019}
(\lambda+m+ar)k'_{r,m}=b_{n_{0}-1}r+\lambda+m.
\end{equation}
Replace $r$ by $-r$, $m$ by $m+r$  in (\ref{5.018}) respectively,
then we have
$$(\lambda+m+r-ar)k_{-r,m+r}=-b_{n_{0}}r+\lambda+m+r.$$
 Similar to the proof of Lemma \ref{l3.5},
\begin{equation}\label{5.020}
(k_{r,m}-1)[(-b_{n_{0}}r+\lambda+m+r)k_{r,m}-(b_{n_{0}}r+\lambda+m)]=0.
\end{equation}
Furthermore,  we can deduce that if $b_{n_{0}}\neq0,1$, then
$k_{r,m}=1$ and $b_{n_{0}}=a$; if $b_{n_{0}}=0$ or $1$, then
$k_{r,m}=1$ or
$k_{r,m}=\dis\frac{b_{n_{0}}r+\lambda+m}{-b_{n_{0}}r+\lambda+m+r}$.

(1)  $\lambda\in\Z$. If $b_{n_{0}}=0$,  let $m=-\lambda$ in
(\ref{5.018}), then we have $ ark_{r,-\lambda}=0$ for all
$r\in\Z$. By the fact that $k_{r, -\lambda}\neq 0$, we have $a=0$.
If $b_{n_{0}}=1$, we have $(a-1)rk_{r,-\lambda-r}=0$ for all
$r\in\Z$  by (\ref{5.018}). Similarly we have $a=1$.
 Therefore, if $b_{n_{0}}=0$ or $1$ and $\lambda\in\Z$, we also
 have $b_{n_{0}}=a$ and $k_{r,m}=1$.

(2)  $\lambda\not\in\Z$. Then $b_{n_{0}}r+\lambda+m\neq0,
-b_{n_{0}}r+\lambda+m+r\neq0$ for all $r,m\in\Z$. By the fact that
$k_{r, m}\neq 0$ for all $r,m\in\Z$ and (\ref{5.018}), we have
$ar+\lambda+m\neq 0$ and
$$k_{r,m}=\frac{b_{n_{0}}r+\lambda+m}{ar+\lambda+m}.$$
Obviously, $b_{n_{0}}=a$ if $k_{r,m}=1$.  On the other hand, if
$k_{r,m}=\dis\frac{b_{n_{0}}r+\lambda+m}{-b_{n_{0}}r+\lambda+m+r}$,
then we have
$$\frac{b_{n_{0}}r+\lambda+m}{-b_{n_{0}}r+\lambda+m+r}=\frac{b_{n_{0}}r+\lambda+m}{ar+\lambda+m}.$$
It is easy to see that $a=1-b_{n_{0}}$  and
\begin{equation}
\begin{cases}\nonumber
 & k_{r,m}=\dis\frac{\lambda+m}{\lambda+m+r},\;\;\; b_{n_{0}}=0;\\
 & k_{r,m}=\dis\frac{\lambda+m+r}{\lambda+m},\;\;\; b_{n_{0}}=1.
\end{cases}
\end{equation}
For $k'_{r,m}$ and $b_{n_{0}-1}$ we have the similar results. So we
have

(i) if $a\neq 0,1$ or  $\lambda\in\mathbb{Z}$, then
$b_{n_{0}}=a=b_{n_{0}-1} \;\; and \; k_{r,m}=k'_{r,m}=1;$

(ii) if $a=0$ or $1$ and $\lambda\not\in\mathbb{Z}$, then
$k_{r,m}=k'_{r,m}$ and

(a) if $k_{r,m}=1$, then $b_{n_{0}}= b_{n_{0}-1}=a$,

(b) if $k_{r,m}=\dis\frac{\lambda+m}{\lambda+m+r}$, then
$b_{n_{0}}=b_{n_{0}-1}=0, a=1$,

(c) if $k_{r,m}=\dis\frac{\lambda+m+r}{\lambda+m}$, then
$b_{n_{0}}=b_{n_{0}-1}=1, a=0$.
\\
Therefore, we always have $k_{r,m}=k'_{r,m}$ and $b_{n}=b$ for all
$r,m,n\in\Z$. Furthermore,
$$g_{r,s}(m,n)=k_{r,m}(as+\mu+n),$$ for all $r,s, m,n\in\Z$.

{\bf Case 1.}\quad   If $k_{r,m}=1$,  then  $b=a$  and
$g_{r,s}(m,n)=as+\mu+n$. Then by (\ref{5.26}), we have
\begin{eqnarray*}
& &f_{r,0}(m,n_{0}-1)=ar+\lambda+m+af(m)+(1-a)f(m+r);
\\
& &f_{r,s}(m,n_{0}-1)=ar+\lambda+m+a(s+1)f(m),\quad\; s\neq0 ;
\\
& &f_{r,s}(m,n_{0}-s-1)=ar+\lambda+m+(1-a)(s+1)f(m+r),\quad\;
s\neq0 ;
\\
& &f_{r,s}(m,n)=ar+\lambda+m, \quad\;n\neq n_{0}-1, n\neq
n_{0}-s-1.
\end{eqnarray*}
Let $h=0, k=-2, s=0, n=n_{0}-1$ in (\ref{e3.4}), then we have
\begin{equation}\label{5.022}
f_{0,-2}(m+r,n_{0})f_{r,0}(m,n_{0}-1)-f_{r,0}(m,n_{0}-2)f_{0,-2}(m,n_{0}-1)=rf_{r,-1}(m,n_{0}-1).
\end{equation}
Then
\begin{equation}\label{5.023}
a[(a+1)r+2\lambda+2m]f(m)+(1-a)(\lambda+m+r)f(m+r)=0.
\end{equation}
Let $r=0$ in (\ref{5.023}),  then we have
\begin{equation}\label{5.024}
(a+1)(\lambda+m)f(m)=0.
\end{equation}
Let $h=0, k=1, s=0, n=n_{0}-1$ in (\ref{e3.4}), and we have
\begin{equation}\label{5.025}
f_{0,1}(m+r,n_{0})f_{r,0}(m,n_{0}-1)-f_{r,0}(m,n_{0}+1)f_{0,1}(m,n_{0}-1)=rf_{r,2}(m,n_{0}-1).
\end{equation}
Then
\begin{equation}\label{5.026}
[(3a-2)(\lambda+m)+(2a-4)ar]f(m)-(a-1)(\lambda+m+r)f(m+r)=0.
\end{equation}
Let $r=0$ in (\ref{5.026}), then
\begin{equation}\label{5.027}
(2a-1)(\lambda+m)f(m)=0.
\end{equation}
By (\ref{5.024}) and (\ref{5.027}), we obtain
\begin{equation}\label{5.028}
(\lambda+m)f(m)=0.
\end{equation}
Then by (\ref{5.023}) and (\ref{5.026}), we know
\begin{equation}\label{5.029}
af(m)=0.
\end{equation}
It is easy to see that if $a\neq 0$ or $\lambda\not\in\Z$,  we
have $f(m)=0$ for all $m\in\Z$. If $\lambda\in\Z$ and $a=0$, by
(\ref{5.028}), we have $f(m)=0$ for all $m\neq-\lambda$. Let
$h=k=-1, r=1, s=0, m=-\lambda, n=n_{0}-1$ in (\ref{e3.4}),  then
we have $f(-\lambda)=0$. Therefore,
$$f(m)=0,$$ for all $m\in\Z$.  We have (\ref{10.1}).

{\bf Case 2.}\quad  If
$k_{r,m}=\displaystyle\frac{\lambda+m+r}{\lambda+m}$, then $b=1$,
$a=0$ and $g_{r,s}(m,n)=\dis\frac{\lambda+m+r}{\lambda+m}(\mu+n)$,
where $\lambda\not\in\Z$.  Note that $\mu+n_{0}=0$. By (\ref{5.26}),
we have
\begin{eqnarray*}
& &f_{r,s}(m,n)=r+\lambda+m,\;\;\; n\neq n_{0}-1, s+n\neq n_{0}-1;
\\
& &f_{r,s}(m,n_{0}-1)=r+\lambda+m,\quad\; s\neq0 ;
\\
& &f_{r,0}(m,n_{0}-1)=r+\lambda+m+
\displaystyle\frac{\lambda+m+r}{\lambda+m}f(m+r);
\\
&
&f_{r,s}(m,n_{0}-s-1)=r+\lambda+m+(s+1)\displaystyle\frac{\lambda+m+r}{\lambda+m}f(m+r),
\quad\;s\neq0.
\end{eqnarray*}
By (\ref{5.022}),  we have
$$(\lambda+m+r)\frac{\lambda+m+r}{\lambda+m}f(m+r)=r(\lambda+m+r).$$
Let $r=0$, since $\lambda\not\in\Z$,  we also have $f(m)=0$ for all
$m\in\Z$. Therefore,  (\ref{10.2}) holds.

{\bf Case 3.}\quad If $k_{r,m}=\dis\frac{\lambda+m}{\lambda+m+r}$,
then $b=0$, $a=1$ and
$g_{r,s}(m,n)=\dis\frac{\lambda+m}{\lambda+m+r}(s+\mu+n)$, where
$\lambda\not\in\Z$.   Note that $\mu+n_{0}=1$.  By (\ref{5.26}), we
have
\begin{eqnarray*}
& &f_{r,s}(m,n)=\lambda+m,\;\;\; n\neq n_{0}-1, s+n\neq n_{0}-1;
\\
& &f_{r,s}(m,n_{0}-s-1)=\lambda+m, \quad\; s\neq0 ;
\\
&
&f_{r,s}(m,n_{0}-1)=\lambda+m+(s+1)\frac{\lambda+m}{\lambda+m+r}f(m),\quad\;
s\neq0 ;
\\
& &f_{r,0}(m,n_{0}-1)=\lambda+m+\frac{\lambda+m}{\lambda+m+r}f(m).
\end{eqnarray*}
By (\ref{5.022}),  we have
$$(\lambda+m+r)[\lambda+m+\frac{\lambda+m}{\lambda+m+r}f(m)]
-(\lambda+m)[\lambda+m-\frac{\lambda+m}{\lambda+m+r}f(m)]=r(\lambda+m).$$
Hence,  $(\lambda+m)f(m)=0.$ Since $\lambda\not\in\Z$, we have
$f(m)=0$ for all $m\in\Z$. Therefore, we have (\ref{10.3}). \qed

\section{\bf Proof of Theorem \ref{t3.1} for the other cases}
\label{Preliminaries}
\label{sub5-c-related}

In this section we prove  that there are no other situations except
the seven ones in Theorem \ref{t3.1}.

\begin{lem}\label{l3.11}
Suppose that
\begin{eqnarray*}
&&Q_{0,s}v_{m,n}=(a_{m}s+\mu+n)v_{m, n+s},\\
&&P_{r,-1}v_{m,n}=(b_{n}r+m)v_{m+r, n}, m\neq 0, m+r\neq 0,\\
&&P_{r,-1}v_{m,n_{0}}=(r+m)v_{m+r, n_{0}}, m\neq 0,\\
&&P_{r,-1}v_{0,n_{0}}=r(1+(r+1)b'_{n_{0}})v_{r, n_{0}},
\end{eqnarray*}
for some $n_{0}\in\mathbb{Z}$. Then $b'_{n_{0}}=0$.
\end{lem}

\pf  By (\ref{5.0}), we have
\begin{equation}\label{11.1}
(\mu+n-a_{m+r})f_{r,-1}(m,n)=(\mu+n-a_{m})f_{r,-1}(m,n-1).
\end{equation}
Let $m\neq0, m+r\neq0$ in (\ref{11.1}), we have
\begin{equation}\label{11.2}
(\mu+n-a_{m+r})(b_{n}r+m)=(\mu+n-a_{m})(b_{n-1}r+m),\;\;\; m\neq0,
m+r\neq0.
\end{equation}
Replace $r$ by $-r$, $m$ by $m+r$ in (\ref{11.2}) respectively,
then we have
\begin{equation}\label{11.3}
(\mu+n-a_{m})(-b_{n}r+m+r)=(\mu+n-a_{m+r})(-b_{n-1}r+m+r),\;\;\;
m\neq0, m+r\neq0.
\end{equation}
From (\ref{11.2}) and (\ref{11.3}), we have
\begin{equation}\label{11.4}
(b_{n}+b_{n-1}-1)(a_{m+r}-a_{m})=0,\;\;for\;all\; m\neq0,
m+r\neq0.
\end{equation}
Note that $b_{n_{0}}=1$. Similar to the proof of Lemma \ref{l3.4} we
have
 $a_{m+r}=a_{m}$ for all $m\neq0,
m+r\neq0$. Therefore,  $$a_{m}=a,\;\;\; m\neq0.$$  By (\ref{11.2}),
we have
\begin{equation}\label{11.7}
(\mu+n-a)(b_{n}-b_{n-1})=0.
\end{equation}
Let $m=0, r\neq 0$ in (\ref{11.1}), then we have
\begin{equation}\label{11.8}
(\mu+n-a)f_{r,-1}(0,n)=(\mu+n-a_{0})f_{r,-1}(0,n-1).
\end{equation}

{\bf Case 1.}\quad If $\mu+n-a\neq 0$ for all $n\in\Z$, then by
(\ref{11.7}), we have $b_{n}=b_{n_{0}}=1$ for all $n\in\Z$. Then
$$f_{r,-1}(m,n)=r+m, \;\;\; m\neq0, m+r\neq0.$$ On the other hand,
for each $n\in\Z$, $V^{R}_{n}=\bigoplus_{m\in\Z}\mathbb{C}v_{m,n}$
is a module of intermediate series over ${\mathcal P}_{-1}$, so we
have
\begin{equation}
\begin{cases}\nonumber
 & f_{r,-1}(m,n)=r+m,\;\;\; m\neq0,
 \\
 & f_{r,-1}(0,n)=r(1+(r+1)b'_{n}).
\end{cases}
\end{equation}
By (\ref{11.8}), we have
$$(\mu+n-a)r(1+(r+1)b'_{n})=(\mu+n-a_{0})r(1+(r+1)b'_{n-1}),$$
for all $r, n\in\Z$. Therefore  $a_{0}=a,  b'_{n}=b'_{n-1},$  and
$$a_{m}=a,\;\;\; b'_{n}=b', $$ for all $m, n\in\Z$.
Similar to the proof of Lemma \ref{l3.5}, we have
$$g_{r,s}(m,n)=k_{r,m}(as+\mu+n),$$
\begin{equation}\label{11.9}
(\mu+n+1-a)f_{0,0}(m,n)=nm+(\mu+1-a)f_{0,0}(m,0),
\end{equation}
\begin{equation}\nonumber
k_{r,m}[(\mu+2-a)f_{0,0}(m,1)-(\mu+1-a)f_{0,0}(m+r,0)+r(\mu+1)]=f_{r,-1}(m,1).
\end{equation}
By (\ref{11.9}), we have
\begin{equation}\label{11.11}
k_{r,m}(\mu+1-a)(f_{0,0}(m,0)-f_{0,0}(m+r,0))+k_{r,m}(m+r+r\mu)=f_{r,-1}(m,1).
\end{equation}
Replace $r$ by $-r$, $m$  by $m+r$  in (\ref{11.11}) respectively,
then we obtain
$$k_{-r,m+r}(\mu+1-a)(f_{0,0}(m+r,0)-f_{0,0}(m,0))+k_{-r,m+r}(m-r\mu)=f_{-r,-1}(m+r,1).$$
According to  $k_{r,m}k_{-r,m+r}=1$, we have
$$k_{r,m}(\mu+1-a)(f_{0,0}(m+r,0)-f_{0,0}(m,0))+k_{r,m}(m-r\mu)=f_{-r,-1}(m+r,1)k_{r,m}^{2}.$$
By (\ref{11.11}), we get
\begin{equation}\label{11.12}
f_{-r,-1}(m+r,1)k_{r,m}^{2}-(2m+r)k_{r,m}+f_{r,-1}(m,1)=0.
\end{equation}
Let $m\neq 0, m+r\neq 0$ in (\ref{11.12}), then
\begin{equation}\nonumber
mk_{r,m}^{2}-(2m+r)k_{r,m}+r+m=0,\;\;\; m\neq 0, m+r\neq 0.
\end{equation}
So
\begin{equation}\label{11.14}
k_{r,m}=1\; or\; k_{r,m}=\frac{r+m}{m},\;\;\; m\neq 0, m+r\neq 0.
\end{equation}
Let $m= 0, r\neq 0$ in (\ref{11.12}), then we obtain
\begin{equation}\label{11.15}
k_{r,0}=1+(r+1)b',\;\;\; r\neq 0.
\end{equation}
By $k_{r,0}k_{-r,r}=1$, we have $k_{-r,r}=\dis\frac{1}{1+(r+1)b'}$
for $r\neq 0$. By the fact that $k_{h+r,m}=k_{h,m+r}k_{r,m}$ for
all $h, r, m\in\Z$, let $m=-r, h\neq 0, r\neq0$, then if
$k_{r,m}=1 $ for $m\neq 0, m+r\neq 0$ we have
$$1=\dis\frac{1}{1+(1-r)b'}[1+(h+1)b'].$$
Then
$$(h+r)b'=0\;\;for\; all\; h\neq 0,\; r\neq 0.$$
Therefore, $b'=0$.  If $k_{r,m}=\dis\frac{r+m}{m}$ for $m\neq 0,
m+r\neq 0$, then we have
$$-\dis\frac{h}{r}=\dis\frac{1}{1+(1-r)b'}(1+(h+1)b'),\;\; h\neq 0,\; r\neq 0.$$
So $b'=-1$,  and
\begin{equation}\nonumber
\begin{cases}
k_{r,m}=\dis\frac{r+m}{m}, & m\neq 0,\; m+r\neq 0,\\
k_{r,0}=-r,          &  r\neq 0,\\
k_{r,-r}=\dis\frac{1}{r}, &   r\neq 0,
\end{cases}
\end{equation}
and
\begin{equation}
\begin{cases}\nonumber
f_{r,-1}(m,n)=r+m,\;\; m\neq 0,\\
f_{r,-1}(0,n)=-r^{2}.
\end{cases}
\end{equation}
Using  (\ref{11.11}), we have
$$f_{0,0}(r,0)-f_{0,0}(0,0)=\dis\frac{\mu}{\mu+1-a}r,$$
$$f_{0,0}(m+r,0)-f_{0,0}(m,0)=\dis\frac{\mu+1}{\mu+1-a}r,\;\;\; m\neq 0,\; m+r\neq 0. $$
Then by (\ref{5.15-}), we obtain
$$f_{r,0}(0,0)=\dis\frac{\mu}{\mu+1-a}r,\;\;\; r\neq 0,$$
$$f_{r,0}(m,0)=\dis\frac{\mu+1}{\mu+1-a}(r+m),\;\;\; m\neq 0,\; r\neq 0. $$
Let $h=r, r=s=k=n=0$ in (\ref{e3.5}), then we have
$$f_{r,0}(m,0)=k_{r,m}[(\mu+1)f_{0,0}(m,0)-\mu
f_{0,0}(m+r,0)+r(a+\mu)].$$ Therefore
\begin{equation}\label{11.16}
f_{0,0}(m,0)=k_{-r,m+r}f_{r,0}(m,0)+\mu[f_{0,0}(m+r,0)-f_{0,0}(m,0)]-r(a+\mu),
\end{equation}
and
\begin{eqnarray*}
&
&f_{0,0}(0,0)=-\dis\frac{\mu}{\mu+1-a}-\frac{\mu-a(a-1)}{\mu+1-a}r,
\\
&
&f_{0,0}(m,0)=\dis\frac{\mu+1}{\mu+1-a}m+\frac{a(a-1)}{\mu+1-a}r,\;\;
m\neq 0.
\end{eqnarray*}
This forces $a=0$ or $a=1$ and $\mu=0$. Thus $\mu+a=0$ or
$\mu-1+a=0$, a contradiction.

{\bf Case 2.}\quad  If there exists $n_{1}\in\Z$ such that
$\mu+n_{1}-a=0$, then by (\ref{11.7}) we have
\begin{equation}
b_{n}=
\begin{cases}\nonumber
 & b_{n_{1}},\;\;\; n\geq n_{1},\\
 & b_{n_{1}-1},\;\;\; n\leq n_{1}-1.
\end{cases}
\end{equation}
 Without loss of generality, we may assume $n_{0}\geq n_{1}$, then
 $b_{n}=1$ for $n\geq n_{1}$.
 By the  results on representations of the Virasoro algebra, we
have
\begin{equation}
\begin{cases}\nonumber
 f_{r,-1}(m, n)=r+m, \;\; m\neq 0,\\
 f_{r,-1}(0, n)=r(1+(r+1)b'_{n}),
\end{cases}\quad\quad\quad\;n\geq n_{1}.
\end{equation}
 In (\ref{11.8}), letting $m=0, r\neq0, n=n_{0}+1=n_{1}+
n'$ , where $ n'> 0$,   we have
$$ n'f_{r,-1}(0,n_{0}+1)=(a-a_{0}+ n')f_{r,-1}(0,n_{0}).$$
Then  $n'r(1+(r+1)b'_{n_{0}+1})=(a-a_{0}+ n')r(1+(r+1)b'_{n_{0}})$
for all  $r\in\Z$. So,
$$b'_{n_{0}+1}=b'_{n_{0}},\;\;\; a_{0}=a.$$ Therefore  $a_{m}=a$
for all $m\in\Z$. By (\ref{11.8}), we have
\begin{equation}
f_{r,-1}(m,n)=
\begin{cases}\nonumber
f_{r,-1}(m,n_{1}),\;\;\;\; n\geq n_{1},\\
f_{r,-1}(m,n_{1}-1),\;\;\; n\leq n_{1}-1,
\end{cases}
\end{equation}
and
\begin{equation}
\begin{cases}\nonumber
f_{r,-1}(m,n)=r+m,\;\;m\neq 0 ,\\
f_{r,-1}(0,n)=r(1+(r+1)b'_{n_{0}}),
\end{cases}\quad\quad\quad\;n\geq n_{1}.
\end{equation}
Similar to the proof of Lemma \ref{l3.10}, we assume
$$g_{r,-1}(m,n)=k_{r,m}(\mu+n-a),\quad for\; n\geq n_{1}.$$
Letting $h=k=0, s=-1, m=-r\neq 0, n\geq n_{1}$ in (\ref{e3.4}), we
have
$$f_{0,0}(0, n)f_{r,-1}(-r, n)-f_{r,-1}(-r, n+1)f_{0,0}(-r, n)=rf_{0,0}(-r, n).$$
Then we can deduce that
$$f_{0,0}(-r, n)=0, \;\;\; r\neq 0, n\geq n_{1}.$$
Letting  $h=-r\neq 0, k=0, s=-1, m=0, n\geq n_{1}$ in
(\ref{e3.4}), we obtain
$$f_{-r,0}(r, n)f_{r,-1}(0, n)-f_{r,-1}(-r, n+1)f_{-r,0}(0, n)=2rf_{0,0}(0, n).$$
Then
$$f_{0,0}(0, n)=0, \;\;\;n\geq n_{1}.$$
Therefore, $$f_{r,0}(-r, n)=0, \;\;\;n\geq n_{1}.$$
 Letting $h=k=0, s=0, m=0, n= n_{1}+1$ in (\ref{e3.5}),
we obtain
\begin{equation}\label{11.16}
f_{0,0}(r,n_{1}+1)g_{r,0}(0,n_{1}+1)+f_{r,0}(0,n_{1}+1)=rg_{r,1}(0,n_{1}+1).
\end{equation}
By (\ref{5.17}), we get
$$g_{r,0}(0, n_{1}+1)=(a+1)k_{r,0},\;\; g_{r,1}(0, n_{1}+1)=(2a+1)k_{r,0},$$
$$g_{-r,-1}(r, n_{1}+2)=2k_{-r,r},\;\; g_{-r,0}(r, n_{1}+1)=(a+1)k_{-r,r}.$$
By (\ref{11.16}), we have
\begin{equation}\label{11.17}
f_{0,0}(r, n_{1}+1)(a+1)k_{r,0}+f_{r,0}(0,
n_{1}+1)=r(2a+1)k_{r,0}.
\end{equation}
Letting $h=k=0, s=-1, m=0, n= n_{1}+1$ in (\ref{e3.4}),  we have
$$
f_{0,0}(r, n_{1}+1)f_{r,-1}(0, n_{1}+1)-f_{r,-1}(0,
n_{1}+2)f_{0,0}(0, n_{1}+1)=rf_{r,0}(0, n_{1}+1).
$$
So
$$r(1+(r+1)b'_{n_{0}})f_{0,0}(r, n_{1}+1)=rf_{r,0}(0, n_{1}+1),$$
and therefore
\begin{equation}\label{11.18}
f_{r,0}(0, n_{1}+1)=(1+(r+1)b'_{n_{0}})f_{0,0}(r, n_{1}+1), \;\;
r\neq 0.
\end{equation}
Letting $h=-r, k=-1, r=0, s=0, m=r\neq 0, n=n_{1}+1$ in
(\ref{e3.5}), we have $$g_{-r,-1}(r,n_{1}+2)f_{0,0}(r,
n_{1}+1)=rg_{-r,0}(r,n_{1}+1).$$
 Then we obtain
\begin{equation}\label{11.19}
f_{0,0}(r, n_{1}+1)=\frac{a+1}{2}r, \;\; r\neq 0.
\end{equation}
So we have $f_{0,0}(m, n_{1}+1)=\dis\frac{a+1}{2}m$ for all
$m\in\Z$, and
$$f_{r,0}(0,n_{1}+1)=\dis\frac{a+1}{2}(1+(r+1)b'_{n_{0}})r,\;\;\;
r\neq 0.$$ On the other hand,
$$f_{0,0}(r,n_{1}+1)(a+1)k_{r,0}+f_{r,0}(0,n_{1}+1)=(2a+1)rk_{r,0}.$$
Then we have
$$(2a+1-a^{2})k_{r,0}=(a+1)(1+(r+1)b'_{n_{0}}),\;\;\; r\neq 0.$$
Letting $h=0, k=-1, r=0, s=0, n=n_{1}+1$ in (\ref{e3.5}), we have
\begin{equation}\label{11.20}
2f_{0,0}(m, n_{1}+1)-f_{0,0}(m, n_{1})=m.
\end{equation}
Then $$f_{0,0}(m, n_{1})=am.$$
 Let $s=0, n=n_{1}+1$ in (\ref{5.16-}), then we get
$$k_{r,m}[2f_{0,0}(m, n_{1}+1)-f_{0,0}(m+r,n_{1})+r(a+1)]=f_{r,-1}(m, n_{1}+1).$$
So
\begin{equation}\label{11.21}
(r+m)k_{r,m}=f_{r,-1}(m,n_{1}+1).
\end{equation}
Let  $m\neq 0$ and $m=0$  in (\ref{11.21}) respectively, then we
have
$$k_{r,m}=1, \;\;\;m\neq 0,\; m+r\neq 0;$$
$$k_{r,0}=1+(r+1)b'_{n_{0}}, \;\;\; r\neq 0.$$
Similar to the proof above, we have $b'_{n_{0}}=0.$ \qed
\\

The proof of the following lemma is similar that of Lemma
\ref{l3.11}.

\begin{lem}\label{l3.13}
Suppose that
\begin{eqnarray*}
&&Q_{0,s}v_{m,n}=(a_{m}s+\mu+n)v_{m, n+s},\\
&&P_{r,-1}v_{m,n}=(b_{n}r+m)v_{m+r, n},\;\;m\neq 0, m\neq -r,\\
&&P_{r,-1}v_{m,n_{0}}=mv_{m+r, n_{0}},\;\; m\neq -r,\\
&&P_{r,-1}v_{-r,n_{0}}=-r(1+(r+1)a')v_{0, n_{0}},
\end{eqnarray*}
for some $n_{0}\in\mathbb{Z}$. Then $a'=0$.
\end{lem}

\begin{lem}\label{l3.22}
Suppose that
\begin{eqnarray*}
&&Q_{0,s}v_{m,n}=(s+n)v_{m, n+s},\;\;\; n\neq 0,\\
&&Q_{0,s}v_{m,0}=s(1+(s+1)a'_{m}))v_{m, s},\\
&&P_{r,-1}v_{m,n}=(b_{n}r+\lambda+m)v_{m+r, n}.
\end{eqnarray*}
Then $a'_{m}=0$ for all $m\in\Z$.
\end{lem}

\pf
By (\ref{5.0}), we have $(n-1)f_{r,-1}(m,n)=(n-1)f_{r,-1}(m,n-1).$
Then $$f_{r,-1}(m,n)=f_{r,-1}(m,n-1),\;\;\; n\neq 1.$$ That is,
\begin{equation}
b_{n}=
\begin{cases}\nonumber
b_{1}, & n\geq 1, \\
b_{0}, & n\leq 0.
\end{cases}
\end{equation}
Letting $h=r, k=-1, r=0, s=-1$ in (\ref{e3.6}),  we have
\begin{equation}\label{22.1}
(n-1)g_{r,-1}(m,n-1)=(n-2)g_{r,-1}(m,n).
\end{equation}
Similar to the proof of Lemma \ref{l3.10}, we can assume
\begin{eqnarray*}
& & g_{r,-1}(m,n)=l_{r,m}(n-1), \;\;\; n\geq1,
\\
& &  g_{r,-1}(m,n)=l'_{r,m}(n-1), \;\;\; n\leq0.
\end{eqnarray*}
 Then we deduce that $l_{r,m}, l'_{r,m}$ have the same properties as
$k_{r,m}, k'_{r,m}$:
\begin{equation}\label{22.2}
l_{r,m}l_{-r,m+r}=1, \;\;\;
l_{h,m+r}l_{r,m}=l_{r,m+h}l_{h,m}=l_{h+r,m},
\end{equation}
\begin{equation}\label{22.3}
l'_{r,m}l'_{-r,m+r}=1,\;\;\;
l'_{h,m+r}l'_{r,m}=l'_{r,m+h}l'_{h,m}=l'_{h+r,m}.
\end{equation}
As the proof of Lemma \ref{l3.10}, we have $l'_{r,m}=l_{r,m}, $ for
all $r,m\in\Z$ and
\begin{equation}\label{22.8}
\begin{cases}
&  g_{r,s}(m,n)=l_{r,m}(s+n), \;\;\; n\neq0;
\\
&  g_{r,s}(m,0)= l_{r,m}s (1+(s+1)a'_{m}).
\end{cases}
\end{equation}
Let  $h=0, k=-1, r=0, s=0$ in (\ref{e3.5}), then we have
$$f_{0,0}(m,n)=\lambda+m, \;\;\; n\neq 0.$$
Assume $f_{0,0}(m,0)=\lambda+m+f(m)$ for all $m\in\Z$. Then
\begin{equation}\label{22.9}
f_{0,0}(m,n)=\lambda+m+\delta_{n,0}f(m).
\end{equation}
Let  $h=r, k=-1, r=0, s=0$ in (\ref{e3.5}), then we obtain
$$nl_{r,m}f_{0,0}(m,n)-(n-1)l_{r,m}f_{0,0}(m+r,n-1)=-rl_{r,m}n+b_{n}r+\lambda+m.
\;\;\; n\neq 0.$$ By (\ref{22.9}), we have
$$nl_{r,m}(\lambda+m+\delta_{n,0}f(m))-(n-1)l_{r,m}(\lambda+m+r+\delta_{n,1}f(m+r))=-rl_{r,m}n+b_{n}r+\lambda+m.$$
Then
\begin{equation}\label{22.10}
(\lambda+m+r)l_{r,m}=b_{n}r+\lambda+m.
\end{equation}
Replacing $r$ by $-r$, $m$ by $m+r$ in (\ref{22.10}) respectively,
we have
$$(\lambda+m)l_{-r,m+r}=-b_{n}r+\lambda+m+r.$$
By (\ref{22.2}), we get
\begin{equation}\label{22.11}
(-b_{n}r+\lambda+m+r)l_{r,m}=\lambda+m.
\end{equation}
Using (\ref{22.10}) and (\ref{22.11}), we obtain
\begin{equation}\label{22.11}
b_{n}l_{r,m}=b_{n}, \;\;\; r\neq 0.
\end{equation}
Then $l_{r,m}=1$  for  $r\neq 0$ or $b_{n}=0$ for all $n\in\Z$. If
$l_{r,m}=1$  for  $r\neq 0$, then $l_{r,m}=1$  for all $r,m\in\Z$.
By (\ref{22.10}), we have $b_{n}=1$ for all $n\in\Z$. If $b_{n}=0$
for all $n\in\Z$, by (\ref{22.10}) we have
$$(\lambda+m+r)l_{r,m}=\lambda+m.$$ It is easy to see that
$\lambda\not\in\Z$  and
$l_{r,m}=\dis\frac{\lambda+m}{\lambda+m+r}$.

{\bf Case 1.} If $l_{r,m}=1$  and $b_{n}=1$ for all $r,m,n\in\Z$,
then
\begin{equation}
\begin{cases}\nonumber
g_{r,s}(m,n)=s+n, \;\;\; n\neq0,
\\
g_{r,s}(m,0)=s (1+(s+1)a'_{m}),
\\
f_{r,-1}(m,n)=r+\lambda+m.
\end{cases}
\end{equation}
By (\ref{5.26}),  we have
\begin{equation}\label{22.12}
f_{r,s}(m,n)=rg_{r,s+1}(m,n)+g_{r,s}(m,n+1)f_{0,0}(m,n)-g_{r,s}(m,n)f_{0,0}(m+r,s+n).
\end{equation}
Then we can deduce that
\begin{eqnarray*}
& & f_{r,s}(m,n)=r+\lambda+m,\;\;\; n\neq 0, -1, -s;
\\
& &
f_{r,s}(m,0)=r+\lambda+m+(s+1)f(m)+(s+1)(2r-s(\lambda+m))a'_{m};
\\
& &f_{r,s}(m,-1)=r+\lambda+m+s(s+1)(\lambda+m)a'_{m};
\\
& &f_{r,s}(m,-s)=r+\lambda+m,\;\;\; n\neq 0,-1.
\end{eqnarray*}
Letting   $h=0, k=-2,  s=0$ in (\ref{e3.4}),  we obtain
\begin{equation}\label{22.13}
f_{0,-2}(m+r,n+1)f_{r,0}(m,n)-f_{r,0}(m,n-1)f_{0,-2}(m,n)=rf_{r,-1}(m,n).
\end{equation}
Let $n=0$  and  $n=1$  in (\ref{22.13}) respectively, then we
obtain
\begin{equation}\label{22.14}
(\lambda+m+r)f(m)+(\lambda+m+r)^{2}a'_{m}=0,
\end{equation}
\begin{equation}\label{22.15}
(\lambda+m)(f(m)+2ra'_{m})=0.
\end{equation}
By (\ref{22.15}), we have $$(\lambda+m)f(m)=0,\;\;\;
(\lambda+m)a'_{m}=0.$$ If $\lambda\not\in\Z$, then $f(m)=0,\;
a'_{m}=0$ for all $m\in\Z$. If $\lambda\in\Z$, then $f(m)=0,\;
a'_{m}=0$ for all $m\neq -\lambda$. Let $m=-\lambda$ in
(\ref{22.14}), then
\begin{equation}\label{22.16}
f(-\lambda)+ra'_{-\lambda}=0,
\end{equation}
for all $r\neq 0$. So $$f(-\lambda)=0,\;\;\; a'_{-\lambda}=0,$$ and
therefore, $a'_{m}=0$ for all $m\in\Z$.

{\bf Case 2.} If $l_{r,m}=\dis\frac{\lambda+m}{\lambda+m+r}$  and
$b_{n}=0$ for all $r,m,n\in\Z$, where $\lambda\not\in\Z$,  then
\begin{equation}
\begin{cases}\nonumber
g_{r,s}(m,n)=\dis\frac{\lambda+m}{\lambda+m+r}(s+n), \;\;\;
n\neq0,
\\
g_{r,s}(m,0)=\dis\frac{\lambda+m}{\lambda+m+r}s (1+(s+1)a'_{m}),
\\
f_{r,-1}(m,n)=\lambda+m.
\end{cases}
\end{equation}
By (\ref{22.12}), we have
\begin{eqnarray*}
& & f_{r,s}(m,n)=\lambda+m,\;\;\; n\neq 0, -1, -s;
\\
& &
f_{r,s}(m,0)=\lambda+m+\frac{s+1}{\lambda+m+r}[f(m)+r(s+2)a'_{m};
\\
&
&f_{r,s}(m,-1)=\lambda+m+\frac{\lambda+m}{\lambda+m+r}s(s+1)a'_{m};
\\
& &f_{r,s}(m,-s)=\lambda+m,\;\;\; n\neq 0,-1.
\end{eqnarray*}
Let $n=1$ in (\ref{22.13}) and note that $\lambda\not\in\Z$, then we
have $2ra'_{m}+f(m)=0$ for all $m, r\in\Z$. Therefore, we have
$$f(m)=0,\; a'_{m}=0,\;\;\; \forall\; m\in\Z.$$
\qed

\begin{lem}\label{l3.28}
Suppose that
\begin{eqnarray*}
&&Q_{0,s}v_{m,n}=(s+n)v_{m, n+s},\;\;\; n\neq 0,\\
&&Q_{0,s}v_{m,0}=s(1+(s+1)a'_{m})v_{m, s},\\
&&P_{r,-1}v_{m,n}=(b_{n}r+m)v_{m+r, n},\;\;\; m\neq 0, m+r\neq 0,\\
&&P_{r,-1}v_{m,n_{0}}=(r+m)v_{m+r, n_{0}},\;\;\; m\neq 0,\\
&&P_{r,-1}v_{0,n_{0}}=r(1+(r+1)b'_{n_{0}})v_{r, n_{0}},
\end{eqnarray*}
for some $n_{0}\in\Z$. Then $b'_{n_{0}}=0$,  $a'_{m}=0$, $b_{n}=1$
for all $m, n\in\Z$.
\end{lem}

\pf  By Lemma \ref{l3.22}, we have
\begin{equation}
\begin{cases}\nonumber
&  g_{r,s}(m,n)=l_{r,m}(s+n), \;\;\; n\neq0,
\\
&  g_{r,s}(m,0)= l_{r,m}s (1+(s+1)a'_{m}),
\end{cases}
\end{equation}
and
\begin{equation}
 f_{r,-1}(m,n)=
\begin{cases}\nonumber
 f_{r,-1}(m,1), & n\geq 1, \\
 f_{r,-1}(m,0), & n\leq 0.
\end{cases}
\end{equation}
Note that $b_{n_{0}}=1$. Without loss of generality, we assume
$n_{0}\geq 1$, then $b_{1}=1$ and
\begin{equation}
\begin{cases}\nonumber
& f_{r,-1}(m,n)=r+m,\;\;\; m\neq 0,\\
& f_{r,-1}(m,n)=r(1+(r+1)b'_{n_{0}}),\\
\end{cases}\quad\quad\quad n\geq1.
\end{equation}
Similar to the proof of Lemma \ref{l3.22}, we have
\begin{equation}\label{28.2}
(m+r)l_{r,m}=f_{r,-1}(m,n).
\end{equation}
Let $n\geq 1$ in (\ref{28.2}), then we obtain
$$(m+r)l_{r,m}=m+r,\;\;\; m\neq 0,$$
$$rl_{r,0}=r(1+(r+1)b'_{n_{0}}).$$
Hence,
$$l_{r,m}=1,\;\;\; m\neq 0,\; m+r\neq 0;$$
$$l_{r,0}=(1+(r+1)b'_{n_{0}}),\;\;\; r\neq 0.$$
Similar to the proof above, we have $b'_{n_{0}}=0$  and
$l_{r,m}=1$ for all $r,m\in\Z$. By (\ref{28.2}), we have $b_{n}=1$
for all $n\in\Z$. Therefore, we can deduce that
$$f_{r,-1}(m,n)=r+m,$$ for all $r, m\in\Z$. By Lemma \ref{l3.22},
we have $a'_{m}=0$ for all $m\in\Z$. \qed

\begin{lem}\label{l3.25}
Suppose that
\begin{eqnarray*}
&&Q_{0,s}v_{m,n}=nv_{m, n+s},\;\;\; n\neq -s,\\
&&Q_{0,s}v_{m,-s}=-s(1+(s+1)a'_{m}))v_{m, 0},\\
&&P_{r,-1}v_{m,n}=(b_{n}r+\lambda+m)v_{m+r, n}.
\end{eqnarray*}
Then $a'_{m}=0$ for all $m\in\Z$. Furthermore, $b_{n}=b$ for all
$n\in\Z$, where $b=0$ or $1$.
\end{lem}

\pf
 By (\ref{5.0}), we have
$$nf_{r,-1}(m,n)=nf_{r,-1}(m,n-1).$$
Then $f_{r,-1}(m,n)=f_{r,-1}(m,n-1)$ for $n\neq 0.$ So we have
\begin{equation}
b_{n}=
\begin{cases}\nonumber
b_{0}, & n\geq 0, \\
b_{-1}, & n\leq -1.
\end{cases}
\end{equation}
Letting  $h=r, k=-1, r=0, s=-1$ in (\ref{e3.6}),  we have
\begin{equation}\label{22.1}
ng_{r,-1}(m,n-1)=(n-1)g_{r,-1}(m,n).
\end{equation}
Similar to the proof of Lemma \ref{l3.22}, assume
\begin{eqnarray*}
& & g_{r,-1}(m,n)=l_{r,m}n, \;\;\; n\geq0,
\\
& &  g_{r,-1}(m,n)=l'_{r,m}n, \;\;\; n\leq-1,
\end{eqnarray*}
then we have $l'_{r,m}=l_{r,m}$ for all $r,m\in\Z$ and
\begin{equation}
\begin{cases}\nonumber
&  g_{r,s}(m,n)=nl_{r,m}, \;\;\; n\neq-s;
\\
&  g_{r,s}(m,-s)=-s(1+(s+1)a'_{m+r})l_{r,m}.
\end{cases}
\end{equation}
Similarly, let  $h=0, k=-1, r=0, s=0$ in (\ref{e3.5}), then we can
deduce
$$f_{0,0}(m,n)=\lambda+m+\delta_{n,0}f(m).$$
Let  $h=r, k=-1, r=0, s=0$ in (\ref{e3.5}), then we have
\begin{equation}\label{25.0}
l_{r,m}[\lambda+m+\delta_{n,0}f(m)-\delta_{n,1}f(m+r)]=b_{n}r+\lambda+m.
\end{equation}
Let $n= 0$ and $n=1$  in (\ref{25.0}) respectively, then we have
$$l_{r,m}(\lambda+m+f(m))=b_{0}r+\lambda+m,$$
$$l_{r,m}(\lambda+m-f(m+r))=b_{0}r+\lambda+m.$$
By the fact that  $l_{r,m}\neq 0$ for all $m,r\in\Z$, we have
$f(m)=-f(m+r)$ for all $m,r\in\Z$. Therefore $$f(m)=0, \;\;\;
\forall\; m\in\Z,$$ and
\begin{equation}\label{25.1}
l_{r,m}(\lambda+m)=b_{0}r+\lambda+m.
\end{equation}
 Let $n\neq 0,1$ in (\ref{25.0}), then we have
$$l_{r,m}(\lambda+m)=b_{n}r+\lambda+m.$$
It is easy to see that $b_{n}=b_{0}=b$ for all $n\in\Z$. Replace $r$
by $-r$, $m$ by $m+r$ in (\ref{25.1}) respectively, then we get
$$l_{-r,m+r}(\lambda+m+r)=-br+\lambda+m+r.$$
By (\ref{22.2}), we have
\begin{equation}\label{25.2}
(-br+\lambda+m+r)l_{r,m}=\lambda+m+r.
\end{equation}
Using  (\ref{25.1}) and (\ref{25.2}), we obtain
\begin{equation}\label{25.3}
l_{r,m}(b-1)=b-1, \;\;\; r\neq 0.
\end{equation}
Then $l_{r,m}=1$  or $b=1$ . If $l_{r,m}=1$  for all $r,m\in\Z$,
then by (\ref{25.1}), we have $b=0$. If $b=1$,   by (\ref{25.1}) we
have $(\lambda+m)l_{r,m}=r+\lambda+m.$  It is easy to deduce that
$\lambda\not\in\Z$  and $l_{r,m}=\dis\frac{r+\lambda+m}{\lambda+m}$
for all $r,m\in\Z$.

{\bf Case 1.} If $l_{r,m}=1$  and $b=0$ , then
\begin{equation}
\begin{cases}\nonumber
g_{r,s}(m,n)=n, \;\;\; n\neq-s,
\\
g_{r,s}(m,-s)=-s(1+(s+1)a'_{m+r}),
\\
f_{r,-1}(m,n)=\lambda+m.
\end{cases}
\end{equation}
By (\ref{5.26}) we have
\begin{equation}\label{25.3}
f_{r,s}(m,n)=rg_{r,s+1}(m,n)+(\lambda+m)g_{r,s}(m,n+1)-(\lambda+m+r)g_{r,s}(m,n).
\end{equation}
Then we can deduce that
\begin{eqnarray*}
& & f_{r,s}(m,n)=\lambda+m,\;\;\; s+n\neq 0,-1;
\\
& &f_{r,s}(m,-s-1)=\lambda+m-(s+1)(r(s+2)+(\lambda+m)s)a'_{m+r};
\\
& &f_{r,s}(m,-s)=\lambda+m+(s+1)(\lambda+m+r)a'_{m+r}.
\end{eqnarray*}
Let $n=2$ in (\ref{22.13}), then  we obtain
\begin{equation}\label{25.4}
(\lambda+m)^{2}a'_{m}=0.
\end{equation}
If $\lambda\not\in\Z$, then $a'_{m}=0$ for all $m\in\Z$. If
$\lambda\in\Z$, then $a'_{m}=0$ for  $m\not\in-\lambda$. Then we
have
\begin{equation}
\begin{cases}\nonumber
& f_{r,s}(-\lambda-r,-s-1)=-r-2(s+1)ra'_{-\lambda},
\\
&  f_{r,s}(m,n)=\lambda+m,\;\;\; {\rm for \ other} \ m,\ n ,
\end{cases}
\end{equation}
and
\begin{eqnarray*}
&&f_{r,s}(m,n)=\lambda+m-\delta_{m+r,-\lambda}\delta_{s+n,-1}2r(s+1)a'_{-\lambda},
\\
&&g_{r,s}(m,n)=n-\delta_{m+r,-\lambda}\delta_{s+n,0}s(s+1)a'_{-\lambda}.
\end{eqnarray*}
Letting $s=1, k=0$ in (\ref{e3.6}), we have
$$(1-2n)\delta_{m+r+h,-\lambda}\delta_{n,-1}a'_{-\lambda}=0.$$
Letting  $n\neq -1, m+h+r\neq -\lambda$, we obtain
$$a'_{-\lambda}=0.$$
Therefore, $a'_{m}=0$ for all $m\in\Z$.

{\bf Case 2.} If $l_{r,m}=\dis\frac{r+\lambda+m}{\lambda+m}$ and
$b=1$, where $\lambda\not\in\Z$, then
\begin{equation}
\begin{cases}\nonumber
g_{r,s}(m,n)=\dis\frac{r+\lambda+m}{\lambda+m}n, \;\;\; n\neq-s,
\\
g_{r,s}(m,-s)=-\dis\frac{r+\lambda+m}{\lambda+m}s(1+(s+1)a'_{m+r}),
\\
f_{r,-1}(m,n)=r+\lambda+m.
\end{cases}
\end{equation}
By (\ref{25.3}), we have
\begin{eqnarray*}
& & f_{r,s}(m,n)=r+\lambda+m,\;\;\; s+n\neq 0, -1;
\\
& &
f_{r,s}(m,-s-1)=r+\lambda+m-\frac{r+\lambda+m}{\lambda+m}(s+1)[s(\lambda+m+r)+2r]a'_{m+r};
\\
&
&f_{r,s}(m,-s)=r+\lambda+m+(s+1)\frac{(r+\lambda+m)^{2}}{\lambda+m}a'_{m+r}.
\end{eqnarray*}
Letting $h=0, k=-2, s=0, n=2$ in (\ref{e3.4}),  we have
$$(\lambda+m+r)a'_{m}=0.$$
Then $a'_{m}=0$ for all $m\in\Z$ by $\lambda\not\in\Z$. \qed

\begin{lem}\label{l3.32}
The following case does not exist:
\begin{eqnarray*}
&&Q_{0,s}v_{m,n}=nv_{m, n+s},\;\;\; n\neq -s,\\
&&Q_{0,s}v_{m,-s}=-s(1+(s+1)a'_{m})v_{m, 0},\\
&&P_{r,-1}v_{m,n}=(b_{n}r+m)v_{m+r, n},\;\;\; m\neq 0, m+r\neq 0,\\
&&P_{r,-1}v_{m,n_{0}}=(r+m)v_{m+r, n_{0}},\;\;\; m\neq 0,\\
&&P_{r,-1}v_{0,n_{0}}=r(1+(r+1)b'_{n_{0}})v_{r, n_{0}},
\end{eqnarray*}
for some $n_{0}\in\Z$.
\end{lem}

\pf  By Lemma \ref{l3.25}, we have
\begin{equation}\label{32.1}
f_{r,-1}(m,n)=
\begin{cases}
& f_{r,-1}(m,0), \;\;\;n\geq0,\\
& f_{r,-1}(m,-1)\;\;\;n\leq-1, \\
\end{cases}
\end{equation}
\begin{equation}\label{32.2}
ml_{r,m}=f_{r,-1}(m,n).
\end{equation}
By (\ref{32.2}), we know $b_{n}=b$ for all $n\in\Z$. Note that
$b_{n_{0}}=1$, so $b=1$ and by the results on representations of the
Virasoro algebra, we have
\begin{equation}
\begin{cases}\nonumber
& f_{r,-1}(m,n)=r+m, \;\;\;m\neq 0,\\
& f_{r,-1}(0,n)=r(1+(r+1)b'_{n_{0}}).
\end{cases}
\end{equation}
Let $m=0$ in (\ref{32.2}), then we have
$$0=f_{r,-1}(0,n)=r(1+(r+1)b'_{n_{0}},$$ a contradiction. \qed

Similarly  we have the following lemmas.

\begin{lem}\label{l3.34}
Suppose that
\begin{eqnarray*}
&&Q_{0,s}v_{m,n}=nv_{m, n+s},\;\;\; n\neq -s,\\
&&Q_{0,s}v_{m,-s}=-s(1+(s+1)a'_{m})v_{m, 0},\\
&&P_{r,-1}v_{m,n}=(b_{n}r+m)v_{m+r, n},\;\;\; m\neq 0, m+r\neq 0,\\
&&P_{r,-1}v_{m,n_{0}}=mv_{m+r, n_{0}},\;\;\; m\neq -r,\\
&&P_{r,-1}v_{-r,n_{0}}=-r(1+(r+1)b'_{n_{0}})v_{0, n_{0}},
\end{eqnarray*}
for some $n_{0}\in\Z$. Then $b'_{n_{0}}=0$, $a'_{m}=0$, $b_{n}=0$
for all $m, n\in\Z$.
\end{lem}

\begin{lem}\label{l3.43}
Suppose that
\begin{eqnarray*}
&&Q_{0,s}v_{m,n}=(a_{m}s+n)v_{m, n+s},\;\;\;n\neq 0, n\neq -s,\\
&&Q_{0,s}v_{m_{0},n}=nv_{m_{0}, n+s},\;\;\;n\neq -s,\\
&&Q_{0,s}v_{m_{0},-s}=-s(1+(s+1)a'_{m_{0}})v_{m_{0}, 0},\\
&&P_{r,-1}v_{m,n}=(b_{n}r+m)v_{m+r, n},\;\;\; m\neq 0, m+r\neq 0,\\
&&P_{r,-1}v_{m,n_{0}}=mv_{m+r, n_{0}},\;\;\; m\neq -r,\\
&&P_{r,-1}v_{-r,n_{0}}=-r(1+(r+1)b'_{n_{0}})v_{0, n_{0}},
\end{eqnarray*}
for some $m_{0}, n_{0}\in\Z$. Then $a'_{m_{0}}=0$, $b'_{n_{0}}=0$
and $a_{m}=0$, $b_{n}=0$ for all $m, n\in\Z$.
\end{lem}

It follows from Lemmas 5.1-5.8 that there are no other situations
except the seven ones in Theorem \ref{t3.1}.




\begin{thebibliography}{9999}
\bibitem{CP} V. Chari, A. Pressley,
{\em  Unitary representations of the Virasoro algebra and a
conjecture of Kac}, Compositio  Mathematica~{\bf67} (1988),
315-342.

\bibitem{DZh} D. Dokovic and K. Zhao, {\em Derivations, isomorphisms and
second cohomology of generalized Witt algebras}, Trans. Amer. Math.
Soc. {\bf 350} (1998), 643-664.


\bibitem{JM} C. Jiang and D. Meng, {\em The derivation algebra of
the associative algebra $\C_{q}[X,Y,X^{-1},Y^{-1}]$}, Communications
in Algebra {\bf 26} (1998), 1723-1736.


\bibitem{K} I. Kaplansky ,
{\em  The Virasoro algebra}, Commun. Math. Phys.~{\bf86} (1982),
149-54.

\bibitem{KS} I. Kaplansky, L. J. Santharoubane,
{\em  Harish Chandra modules over the Virasoro algebra},
Publ.Math.Sci.Res.Inst.~{\bf4} (1987), 217-231.

\bibitem{KPS} E. Kirkman, C. Procesi and L. Small,
{\em A q-analog for the Virasoro algebra}, Communications in
Algebra ~{\bf  22} (1994),  3755-3774.



\bibitem{LT} W. Lin and S. Tan, {\em Non-zero level Harish-Chandra
modules over the Virasoro-like algebra}, J. Pure Appl. Algebra
{\bf 204} (2006), 90-105.



\bibitem{MP} C. Martin, A. Piard,
{\em  Indecomposable modules over the Virasoro Lie algebra and a
Conjecture of V.Kac}, Commun. Math. Phys.~{\bf137} (1991),
109-132.

\bibitem{M} O. Mathieu ,
{\em  Classification of Harish-Chandra modules over the Virasoro
Lie algebra}, Invent. Math.~{\bf107} (1992), 225-234.



\bibitem{PZ} J. Patera, H. Zassenhaus, {\em The higher rank Virasoro algebras,\/}
  Commun. Math. Phys. {\bf 136}(1991), 1-14.

\bibitem{Su1} Y. Su, {\em Generalized Virasoro and super-Virasoro algebras
  and modules of the intermediate series\/}, J. Algebra {\bf 252}(2002),
1-19.


\bibitem{Su2} Y. Su, {\em Simple modules over the high rank Virasoro algebras\/} ,
  Comm. Alg. {\bf 29}(2001),
2067-2080.

\bibitem{SXZ}  Y. Su, X. Xu,  H. Zhang, {\em Derivation-simple algebras and the
   structures of Lie algebras of Witt type\/},  J. Algebra  {\bf  233}(2000),
642-662.



\bibitem{SZ2}   Y. Su, K. Zhao, {\em Second cohomology group of generalized Witt type
  Lie algebras and certain representations\/},  Comm. Alg.  {\bf 30}(2002),
3285-3309.

\bibitem{SZ}  Y. Su, J. Zhou, {\em Some representations of non-graded Lie algebras of
  generalized Witt type\/},   J. Algebra {\bf 246}(2001),
721-738.

\bibitem{X} X. Xu, {\em New generalized simple Lie algebras of Cartan type over a field
  with characteristic 0\/},   J. Algebra {\bf 224}(2000),
 23-58.

\bibitem{Zh} K. Zhao, {\em Weight modules over generalized Witt
algebras with 1-dimensional weight spaces}, Forum Math. {\bf 16}
(2004), 725-748.



\end{thebibliography}
\end{document}